\documentclass{conm-p-l}
\usepackage{amssymb, eucal, eufrak}

\newtheorem{theo}{Theorem}

\newtheorem{prop}[theo]{Proposition}
\newtheorem{defi}{Definition}
\newcommand{\Der}{{\rm Der\:}}
\newcommand{\ad}{{\rm ad\:}}
\newcommand{\gr}{{\rm gr\:}}
\newcommand{\Z}{\mathbb Z}
\newcommand{\rad}{{\rm rad\:}}

\begin{document}

\title[Classification of finite dimensional simple Lie algebras]
{Classification of finite dimensional simple Lie algebras\\ in
prime characteristics}

\author{Alexander Premet}
\address{School of Mathematics, The University of Manchester,
Oxford Road, M13 9PL, United Kingdom}
\email{sashap@maths.man.ac.uk}
\thanks{}

\author{Helmut Strade}
\address{Fachbereich Mathematik, Universit{\"a}t Hamburg, Bundesstrasse 55,
20146 Hamburg, Germany } \email{strade@math.uni-hamburg.de}

\thanks{}


\date{}

%\dedicatory{Dedicated to the memory of Alexei Ivanovich Kostrikin}


\begin{abstract}
We give a comprehensive survey of the theory of finite dimensional
Lie algebras over an algebraically closed field of prime
characteristic and announce that the classification of all finite
dimensional simple Lie algebras over an algebraically closed field
of characteristic $p>3$ is now complete. Any such Lie algebra is
up to isomorphism either classical or a filtered Lie algebra of
Cartan type or a Melikian algebra of characteristic $5$.
\end{abstract}
\maketitle

Unless otherwise specified, all Lie algebras in this survey are
assumed to be finite dimensional. In the first two sections, we
review some basics of modular Lie theory including absolute toral
rank, generalized Winter exponentials, sandwich elements, and
standard filtrations. In Section 3, we give a systematic
description of all known simple Lie algebras of characteristic
$p>3$ with emphasis on graded and filtered Cartan type Lie
algebras. We also discuss the Melikian algebras of characteristic
$5$ and their analogues in characteristics $3$ and $2$. Our main
result (Theorem 7) is stated in Section 4 which also contains
formulations of several important theorems frequently used in the
course of classifying simple Lie algebras. The main principles of
our proof of Theorem 7, with emphasis on the rank two case, are
outlined in Section 5. As suggested by the referee, we mention in
Section 6 some interesting open problems related to the subject.

We would like to thank the referee for careful reading and
valuable comments.

\section{The beginnings}

The theory of Lie algebras over a field $F$ of characteristic
$p>0$ was initiated by Jacobson, Witt and Zassenhaus. In
\cite{J1}, Jacobson investigated purely inseparable field
extensions $E/F$ of the form $E=F(c_1,\ldots,c_n) $ where
$c_i^p\in F$ for all $i\le n$. Although such field extensions do
not possess nontrivial $F$-automorphisms, Jacobson developed for
them a version of Galois theory. The r{\^o}le of Galois
automorphisms in his theory was played by $F$-derivations.

The set ${\rm Der}_F\,E$ of all $F$-derivations of $E$ carries the
following three structures:
\begin{itemize}
\item a natural structure of a vector space over $E$, \item a
natural $p$-structure given by the $p$th power map $D \mapsto
D^p$, \item a Lie algebra structure given by the commutator
product.
\end{itemize}
Let $\mathfrak F$ denote the set of all subfields of $E$
containing $F$ and $\mathfrak L$ the set of all $E$-subspaces of
${\rm Der}_{F}\, E$ stable under the $p$th power map and Lie
bracket in ${\rm Der}_{F}\, E$. Both sets $\mathfrak F$ and
$\mathfrak L$ are partially ordered by inclusion. Given a subset
$X$ in ${\rm Der}_F\,E$ we let $E^X$ denote the subfield of $E$
consisting of all $\alpha\in E$ satisfying $x(\alpha)=0$ for all
$x\in X$.
\begin{theo}[\cite{J1}] The map $\mathfrak{L}\ni L\mapsto E^L\in {\mathfrak F}$ is an
order-reversing bijection between $\mathfrak L$ and $\mathfrak F$.
\end{theo}
\noindent Jacobson singled out the $p$-structure above as being of
major importance for Lie theory.
\begin{defi}[\cite{J1}] A Lie algebra $L$ over $F$ is
called {\em restrictable} if
for any $x\in L$ the derivation $({\rm ad}\, x)^p$ of $L$ is inner.
\end{defi}
\noindent Any restrictable Lie algebra $L$ carries a $p$-mapping
$x\mapsto x^{[p]}$ which enjoys the three following properties:
\begin{enumerate}
\item[1.\ \,] $(\lambda x)^{[p]}=\lambda^p x^{[p]}$, \item[2.\ \,]
$({\rm ad}\, x)^p={\rm ad}\, x^{[p]}$, \item[3.\ \,]
$(x+y)^{[p]}=x^{[p]}+y^{[p]}+\sum_{i=1}^{p-1} s_{i}(x,y),$ where
$s_{i}(x,y)\in L$ are such that
$$\sum_{i=1}^{p-1}is_{i}(x,y)t^{i-1}\,=\,({\rm ad}\,
(tx+y))^{p-1}(x)$$
\end{enumerate}
(here $x,y\in L$, $\lambda\in F$, and $t$ is a variable). Such a
$p$-mapping is uniquely determined up to a $p$-linear map from $L$
into its center $\mathfrak{z}(L)$. It is therefore unique for any
restrictable Lie algebra $L$ with $\mathfrak{z}(L)=(0)$. Once the
mapping $[p]$ is fixed, the pair $(L, [p])$ is called a {\em
restricted Lie algebra}. If $I$ is a {\it restricted} ideal of
$L$, that is an ideal of $L$ such that $I^{[p]}\subset I$, then
the quotient Lie algebra $L/I$ carries a {\it natural} $p$-mapping
given by $(x+I)^{[p]}=x^{[p]}+I$ for all $x\in L$. We mention for
completeness that the Lie algebras of linear algebraic groups over
$F$ are all equipped with canonical $p$-mappings, hence carry {\it
canonical} restricted Lie algebra structures.

From now on we assume that $F$ is algebraically closed. Some time
before 1939  Witt discovered (for any $p>3$) a $p$-dimensional
simple Lie algebra with no finite dimensional analogues in
characteristic 0. The {\em Witt algebra} $W(1;\underline{1})$ has
basis $\{e_{-1},e_0, e_{1},\ldots, e_{p-2}\}$ over $F$ and the Lie
product in $W(1;\underline{1})$ is given by
$$[e_i,e_j]=\left\{\begin{array}{ll}(j-i)e_{i+j}& \mbox{if }\
-1\leq i+j\leq p-2,\\ 0 &\mbox{otherwise}.
\end{array}\right.$$
As Witt himself never published his example, we have only indirect
information about his discovery. Zassenhaus generalized Witt's
example by considering a  subgroup $G$ of order $p^n$ in the
additive group of $F$ and by giving a $p^n$-dimensional vector
space $W_G:=\bigoplus_{g\in G}\, Fe_g\,$ a Lie algebra structure
via $[e_g,e_h] := (h-g)e_{g+h}$ for all $g, h \in G$. Such Lie
algebras are often referred to as {\em Zassenhaus algebras}.

In \cite{Z}, Zassenhaus investigated irreducible representations
of nilpotent Lie algebras over fields of prime characteristics.
This paper is the starting point of the modular representation
theory of Lie algebras.

In \cite{C}, Chang described all irreducible representations of
the Witt algebra $W(1;\underline{1})$. According to \cite{C}, Witt
used the following realization of the Lie algebra
$W(1;\underline{1})$: Let ${\mathcal O}(1;\underline{1})$ denote
the truncated polynomial algebra $F[X]/(X^p)$, and let $x$ be the
image of $X$ in ${\mathcal O}(1;\underline{1})$. Give ${\mathcal
O}(1;\underline{1})$ an algebra structure by setting $\{f,g\} :=
f(dg/dx) - g(df/dx)$ for all $f, g\in {\mathcal
O}(1;\underline{1})$. It is readily seen that the map $e_i\mapsto
x^{i+1}$ extends to an algebra isomorphism
$W(1;\underline{1})\stackrel{\sim}{\longrightarrow} \big({\mathcal
O}(1;\underline{1}),\{\,\cdot\,,\,\cdot\,\}\big)$. For $i\in
\mathbb{F}_p$ set $u_i=(1+x)^{i+1}$. Then $\{u_i,u_j\}
=(j-i)u_{i+j}$ for all $i, j\in\mathbb{F}_p$. This shows that
$W(1;\underline{1})$ is isomorphic to the Zassenhaus algebra
associated with the additive subgroup $\mathbb{F}_p\subset F$.

\section{Some basics}

This section is a short introduction into the general theory of
modular Lie algebras with emphasis on results and techniques used
in Classification Theory. Most of the results discussed here are
valid for any prime $p$.

\subsection{Maximal tori in restricted Lie algebras}

Let $\mathfrak g$ be a restricted Lie algebra over $F$. An element
$x\in \mathfrak{g}$ is called {\it semisimple} (respectively, {\it
nilpotent}) if $x$ lies in the restricted subalgebra of
$\mathfrak{g}$ generated by $x^{[p]}$ (respectively, if
$x^{{[p]}^e}=0$ for $e\gg 0$). For any $x\in\mathfrak{g}$ there
exist unique commuting $x_{s}$ and $x_{n}$ in $\mathfrak g$ such
that $x_s$ is semisimple, $x_n$ is nilpotent, and $x=x_{s}+x_n$.
We denote by $\mathfrak{g}_{x}^0$ the set of all
$y\in\mathfrak{g}$ such that $(\mbox{ad}\,
x)^{\dim\mathfrak{g}}(y)=0$, and define ${\rm
rk}(\mathfrak{g})\,:=\,\min\,\{\dim\mathfrak{g}_{x}^0\,\vert\,x\in
\mathfrak{g}\}.$ If $\dim\,\mathfrak{g}_{x}^{0}\,=\,
\mbox{rk}(\mathfrak{g})$ then $\mathfrak{g}_{x}^0$ is a Cartan
subalgebra of $\mathfrak g$ (this is a standard fact of Lie
theory).

An element $t\in\mathfrak g$ is called ${\it toral}$ if
$\,t^{[p]}=t$. A restricted subalgebra $\mathfrak t$ of $\mathfrak
g$ is called {\it toral} (or a {\it torus} of $\mathfrak g$) if
the $p$-mapping is invertible on $\mathfrak t$. Any toral
subalgebra of $\mathfrak g$ is abelian and admits a basis
consisting of toral elements. Set
$$MT(\mathfrak{g})\,:=\,\max\,\{\dim \mathfrak{t}\,\vert\,\mathfrak{t}\
\mbox{is a torus in }\,\mathfrak{g}\}.$$ A torus $\mathfrak{t}$ of
$\mathfrak g$ is called {\it maximal} if the inclusion $
\mathfrak{t}\subset\mathfrak{t}'$ with $\mathfrak{t}'$ toral
implies $\mathfrak{t}=\mathfrak{t}'$. The centralizer $
\mathfrak{c}_{\mathfrak g}(\mathfrak{t})$ of any maximal torus in
$\mathfrak g$ is a Cartan subalgebra of $\mathfrak g$ and,
conversely, the semisimple elements of any Cartan subalgebra of
$\mathfrak g$ lie in its center and form a maximal torus in
$\mathfrak g$. The reader should be warned, however, that maximal
tori (and their centralizers) in a restricted Lie algebra may have
different dimensions (see \cite{St1}). In other words, there may
exist maximal tori in $\mathfrak g$ of dimension less that
$MT(\mathfrak{g})$.

Let $\mathfrak t$ be a maximal torus of $\mathfrak g$,
$\mathfrak{h}=\mathfrak{c}_{\mathfrak g}(\mathfrak{t})$, and let
$V$ be a finite dimensional restricted $\mathfrak g$-module (this
means that $\rho_V(x^{[p]})=\rho_V(x)^p$ for any $x\in\mathfrak g$
where $\rho_V$ denotes the corresponding representation). Since
$\rho_V(\mathfrak t)$ is abelian and consists of semisimple
elements, $V$ decomposes into weight spaces relative to $\mathfrak
t$:
$$V\,=\,\bigoplus_{\lambda\in{\mathfrak{t}}^*}\,V_{\lambda},\qquad \
\quad V_{\lambda}\,=\, \{v\in V\,\vert\,t\,.\, v=\lambda(t)v\ \
\,\forall\,t\in\mathfrak t\}.$$ The set of $\mathfrak t$-weights
$\{\lambda\in {\mathfrak t}^*\,\vert\,V_{\lambda}\ne 0\}$ of $V$
will be denoted by $\Gamma^{w}(V,\mathfrak{t})$. It is worth
mentioning that if $t$ is a toral element of $\mathfrak t$ then
$\lambda(t)\in\mathbb{F}_p$ for any $\lambda\in\Gamma^w(V,
\mathfrak{t})$. Set $\Gamma(V,\mathfrak{t})=\Gamma^{w}(V,
\mathfrak{t})\setminus\{0\}$. For $\mathbb{F}_p$-independent
linear functions $\mu_1,\ldots,\mu_k\in\Gamma(V,\mathfrak{t})$
define
$$V(\mu_1,\ldots,\mu_k)\,:=\,\bigoplus_{(i_{1},\ldots,
i_k)\in\mathbb{F}_{p}^k}\,V_{i_1\mu_1+\cdots+i_k\mu_k}.$$ The
subspace $V(\mu_1,\ldots,\mu_k)$ is called a $k$-{\it section} of
$V$.

If $V$ is an algebra over $F$ (not necessarily associative or Lie)
and $\mathfrak g$ acts on $V$ as derivations then
$V(\mu_1,\ldots,\mu_k)$ is a subalgebra of $V$. If $V=\mathfrak
g$, the adjoint $\mathfrak g$-module, then $\Gamma=\Gamma(
\mathfrak{g},\mathfrak{t})$ is nothing but the set of roots of
$\mathfrak g$ relative to $\mathfrak t$, and
$$\mathfrak{g}\,=\,\mathfrak{h}\,\oplus\,\sum_{\alpha\in
\Gamma}\,\mathfrak{g}_{\alpha}$$ is the root space decomposition.
A Cartan subalgebra $\mathfrak h$ of $\mathfrak g$ is called {\it
regular} if $\mathfrak{h}\,=\,\mathfrak{c}_{\mathfrak
g}(\mathfrak{t})$ where $\mathfrak t$ is a torus of  maximal
dimension in $\mathfrak g$. By the main result of \cite{P86b}, all
regular Cartan subalgebras of $\mathfrak g$ have dimension equal
to ${\rm rk}(\mathfrak{g})$.

Let $\mathfrak{h}=\mathfrak{c}_{\mathfrak g}(\mathfrak{t})$ be a
regular Cartan subalgebra of $\mathfrak g$. In \cite{Wi}, Winter
proved that for any $x\in\mathfrak{g}_\gamma$ satisfying
$x^{[p]}=0$ the exponential operator $\exp{\rm ad}\,x\in{\rm
GL}(\mathfrak{g})$ maps the root space decomposition of $\mathfrak
g$ relative to $\mathfrak h$ onto that of another regular Cartan
subalgebra, denoted $\mathfrak{h}_x$. To appreciate this result
one should keep in mind that in characteristic $p$ the condition
$x^{[p]}=0$ does not always guarantee that $\exp{\rm ad}\,x$ is an
automorphism of $\mathfrak g$ (for example, consider the case
where $\mathfrak{g}=W(1;\underline{1})$ and $x=e_{-1}$).

In \cite{Wil83}, Wilson assigned a generalized exponential
operator to any root vector $x\in\mathfrak{g}_\gamma$ such that
$x^{[p]}\in\mathfrak{t}$. Inspired by Wilson's construction, the
first author assigned generalized exponential operators to all
root vectors in $\mathfrak{g}_\gamma$; see \cite{P86b}.
Generalized exponential operators and resulting switchings of
regular Cartan subalgebras in $\mathfrak g$ play an important
r{\^o}le in Classification Theory.

Let $\xi\in\mbox{Hom}_{\mathbb{F}_p}(F,F)$ be such that
$\xi^p-\xi={\rm Id}_F$. As $F$ is algebraically closed, it is
straightforward to see that $\xi:F\rightarrow F$ exists and is
uniquely determined up to a linear map from $F$ to $\mathbb{F}_p$.
Given $x\in\mathfrak{g}_\gamma$, where $\gamma\in\Gamma$, we
denote by $m=m(x)$ the least positive integer $k$ with
$x^{[p]^k}\in\mathfrak{t}$  (such an integer exists because
$\mathfrak t$ is a maximal torus in $\mathfrak g$). Set
\[
q(x)=\left\{
\begin{array}{ll}
\sum_{i=1}^{m-1}x^{{[p]}^i}& \mbox{for \   $m>1$}, \\ 0 &
\mbox{for \ $m=1$}.
\end{array}
\right.
\]
Note that $q(x)\in \mathfrak{h}$. Define the {\it generalized
Winter exponential} $E_{x,\xi}\in\mbox{GL}(\mathfrak{g})$ by
setting
$$E_{x,\xi}(y)\,=\,-\sum_{i=0}^{p-1}\,
\prod_{j=i+1}^{p-1}\Big((\xi(\alpha(x^{[p]^m}))+j)\,\mbox{Id}_{\mathfrak
g}-\mbox{ad}\,\, q(x)\Big)\,(\mbox{ad}\,\,x)^i(y)$$ for all $y\in
\mathfrak{g}_\alpha$, where $\alpha\in\Gamma\cup\{0\}$, and
extending to $\mathfrak g$ by linearity (our convention here is
that $\mathfrak{g}_0=\mathfrak{h}$). Notice that if $x^{[p]}=0$
then $E_{x,\xi}=\exp \mbox{ad}\,x$. In general, $E_{x,\xi}$ is a
polynomial in $\ad x$; see \cite{P89}.

According to \cite{P86b}, $\mathfrak{h}_x\,=\,E_{x,\xi}
(\mathfrak{h})$ is a regular Cartan subalgebra of $\mathfrak g$
and
$$\mathfrak{g}\,=\,\mathfrak{h}_x\oplus\,\sum_{\alpha\in\,\Gamma}\,E_{x,\xi}(
\mathfrak{g}_\alpha)$$ is the root space decomposition of
$\mathfrak g$ relative to $\mathfrak{h}_x$. For $t\in\mathfrak{t}$
set
$$t_x:=\,t-\gamma(t)(x+q(x)).$$ The subspace $
\mathfrak{t}_x\,=\,\{t_x\,\vert\,t\in\mathfrak{t}\}$ coincides
with the unique maximal torus in $\mathfrak{h}_x$; see
\cite{P86b}. The set of roots $\Gamma(\mathfrak{g},
\mathfrak{t}_x)$ of $\mathfrak g$ relative to $\mathfrak{t}_x$ has
the form $\Gamma(\mathfrak{g},\mathfrak{t}_x)
\,=\,\{\alpha_{x,\xi}\,\vert\,\alpha\in\Gamma\} \subset
\mathfrak{t}_{x}^*$ where
$$\alpha_{x,\xi}(t_x)\,=\,\alpha(t)-\xi(\alpha(x^{[p]^m}))\,\gamma(t)\qquad\
(\forall\,t_x\in\mathfrak{t}_x).$$ If
$\mathfrak{h}'=E_{y,\xi}(\mathfrak{h})$ for some
$y\in\bigcup_{\,\gamma\in\Gamma}\,\mathfrak{g}_\gamma$, we say
that $\mathfrak{h}'$ is obtained from $\mathfrak h$ by an {\it
elementary switching}. By \cite{P89}, any two regular Cartan
subalgebras of $\mathfrak g$ can be obtained each from another by
a finite sequence of elementary switchings. This result has the
following important consequence:

\begin{prop}[\cite{PS2}] Let $\mathfrak{t}_1$ and $\mathfrak{t}_2$
be two tori of maximal dimension in $\mathfrak g$, $V$ a finite
dimensional restricted $\mathfrak g$-module,
$\Delta_{i}=\Gamma^w(V,\mathfrak{t}_i)$, and $Q_i$ the $
\mathbb{F}_p$-span of $\Delta_i$ in $\mathfrak{t}^*_i$, where
$i=1,2$. Then there exists a linear isomorphism of $
\mathbb{F}_p$-spaces $\psi\,\colon\, Q_1\rightarrow Q_2$ such that
$\psi(\Delta_1)=\Delta_2$ and $\dim\,V_{\mu}=\dim\,V_{\psi(\mu)}$
for all $\mu\in \Delta_1$.
\end{prop}
As a consequence one obtains that
$\mathbb{F}_p^*\,\delta_1\subset\Delta_1$ for some $\delta_1\in
\mathfrak{t}^*_1$ if and only if $
\mathbb{F}_p^*\,\delta_2\subset\Delta_2$ for some $\delta_2\in
\mathfrak{t}^*_2$. Also, $0\in\Delta_1$ if and only if $0
\in\Delta_2$.

\subsection{Absolute toral rank}

It is often useful to view a Lie algebra as a subalgebra of a
restricted Lie algebra.
\begin{defi} [\cite{St-F}] Let $L$ be a Lie algebra. A triple
$({\mathcal L},[p],i)$ where ${\mathcal L}$ is a restricted Lie
algebra with $p$-mapping $[p]\,:\, {\mathcal L}\rightarrow {\mathcal
L}$ and $i\,\colon\, L\hookrightarrow \mathcal L$ is an injective
Lie algebra homomorphism, is called a {\rm $p$-envelope} of $L$ if
the restricted Lie subalgebra of $\mathcal L$ generated by $i(L)$
coincides with $\mathcal L$.
\end{defi}

\noindent The Lie algebra $L$ is often identified with $i(L)\subset
\mathcal L$. We list below a few basic properties of $p$-envelopes.
All proofs can be found in \cite{St-F, St04}.
\begin{itemize}
\item[{\bf 2.2.1.}] Let $({\mathcal L},[p],i)$ and $({\mathcal
L}', [p]',i')$ be two $p$-envelopes of $L$. Then there exists an
isomorphism of restricted Lie algebras $\psi\,\colon\,{\mathcal
L}/\mathfrak{z}({\mathcal
L})\stackrel{\sim}{\longrightarrow}{\mathcal L}'/
\mathfrak{z}({\mathcal L}')$ such that $\psi\circ\pi\circ
i=\pi'\circ i'$ where $\pi$ and $\pi'$ denote the canonical
homomorphisms of restricted Lie algebras ${\mathcal
L}\twoheadrightarrow {\mathcal L}/{\mathfrak z}({\mathcal L})$ and
${\mathcal L}'\twoheadrightarrow {\mathcal
L}'/\mathfrak{z}({\mathcal L}')$. \item[{\bf 2.2.2.}] A
$p$-envelope $({\mathcal L}, [p],i)$ of $L$ is called {\it
minimal} if $\mathfrak{z}({\mathcal L})$ is contained in
$\mathfrak{z}(i(L))$. Any $L$ admits a minimal $p$-envelope, and
any two minimal $p$-envelopes of $L$ are isomorphic as {\it
ordinary} Lie algebras. \item[{\bf 2.2.3.}] Suppose $L$ is
semisimple. Then $L$ has one ``obvious'' minimal $p$-envelope,
namely, the restricted Lie subalgebra of $\Der L$ generated by
${\rm ad}\, L$. This $p$-envelope is semisimple. Any two
semisimple $p$-envelopes of $L$ are isomorphic as {\it restricted}
Lie algebras.
\end{itemize}
\begin{defi}
Let $({\mathcal L},[p],i)$ be a $p$-envelope of $L$. The {\rm
absolute toral rank} of $L$, denoted $TR(L)$, is the maximal
dimension of tori in the restricted Lie algebra ${\mathcal
L}/\mathfrak{z}({\mathcal L})$. In other words,
$$TR(L)\,:=\,MT({\mathcal L}/\mathfrak{z}({\mathcal L})).$$
\end{defi}
In view of ({\bf 2.2.1}), this definition is independent of the
choice of a $p$-envelope of $L$. For $L$ semisimple,
$TR(L)\,=\,MT(L_p)$ where $L_p$ stands for the restricted Lie
subalgebra of $\Der L$ generated by ${\rm ad}\,L$ (see ({\bf
2.2.3})). We shall need a few basic properties of $TR(L)$ all of
which can be found in \cite{St04}:
\begin{itemize}
\item[{\bf 2.2.4.}] $L$ is nilpotent if and only if $TR(L)=0$.
\item[{\bf 2.2.5.}]  If $I$ is an ideal of $L$ then
$TR(L/I)+TR(I)\le TR(L)$. \item[{\bf 2.2.6.}] Let $T$ be a torus
of maximal dimension in a finite dimensional $p$-envelope of $L$
and let $\gamma_1,\ldots,\gamma_k$ be $\mathbb{F}_p$-independent
roots in $\Gamma(L,T)$. Then $$TR(L(\gamma_1,\ldots,\gamma_k))\le
k.$$ In particular, $TR(L(\alpha)) \le 1$ for any
$\alpha\in\Gamma(L,T)$ and $TR(L(\alpha,\beta))\le 2$ for any two
$\alpha,\beta\in\Gamma(L,T)$.
\end{itemize}

\subsection{Sandwich elements}

Given an arbitrary Lie algebra $L$ over a field we define ${\mathcal
S}(L):=\{s\in L\,\vert\, ({\rm ad}\,s)^2=0\}$. The set ${\mathcal
S}(L)$ plays a crucial r{\^o}le in Kostrikin's work on the
restricted Burnside problem (see \cite{Ko90}). If $2L=L$ and
$s\in{\mathcal S}(L)$, then $({\rm ad}\, s)({\rm ad}\, x)({\rm ad}\,
s)=0$ for any $x\in L$. Because of this property the elements of
${\mathcal S}(L)$ are often referred to as {\it sandwich} elements
(the term is due to Kostrikin). As an example, ${\mathcal
S}(W(1;\underline{1}))=\bigoplus_{2i>p}\,Fe_i.$ In general, $S(L)$
is not closed under vector addition however. If $2L=L$, then
${\mathcal S}(L)$ is closed under Lie multiplication (see
\cite{Ko90} for more detail).

Assume until the end of this subsection that $\mbox{char}\, F=p>2$
and let $L$ be finite dimensional over $F$. Let $c\in{\mathcal
S}(L)$ and $x\in L$. Since
$(\mbox{ad}\,c)(\mbox{ad}\,x)(\mbox{ad}\,c)=0$ we have
$(({\mbox{ad}\,c})(\mbox{ad}\,x))^2=0$. This implies that
$\mbox{tr}\,(\mbox{ad}\,c)(\mbox{ad}\,x)=0$. As a consequence,
${\mathcal S}(L)$ is contained in the radical of the Killing form of
$L$. The Lie algebras over $F$ containing nonzero sandwich elements
are called {\it strongly degenerate} (the term is due to Kostrikin).
It follows from the preceding remark that the Killing form of any
strongly degenerate simple Lie algebra over $F$ is identically zero.

By the Engel--Jacobson theorem, the linear span $\langle{\mathcal
S}\rangle$ of ${\mathcal S}={\mathcal S}(L)$ is a nilpotent Lie
subalgebra of $L$. Since ${\mathcal S}(L)$ is invariant under all
automorphisms of $L$ the same is true for the normalizer of
$\langle\mathcal S\rangle$ in $L$. As a consequence, every strongly
degenerate simple Lie algebra $L$ contains a proper nonzero
subalgebra invariant under all automorphisms of $L$. (This remark
also shows that in characteristic $0$ the equality ${\mathcal
S}(L)=\{0\}$ is equivalent to the semisimplicity of $L$.) For $p>3$,
the Lie algebras $L$ over $F$ with ${\mathcal S}(L)=\{0\}$ are
closely related to the Lie algebras of semisimple algebraic groups
over $F$; see the discussion in (3.1) for more detail.

In \cite{KS66}, Kostrikin and Shafarevich conjectured that for
$p>5$ the normaliser of $\langle {\mathcal S}\rangle$ in any
strongly degenerate simple Lie algebra $L$ is a {\it maximal}
subalgebra of $L$. In his PhD thesis and a subsequent series of
preprints, S.A.~Kirillov verified this conjecture for all known
finite dimensional simple Lie algebras of characteristic $p>3$.
Unfortunately, all attempts to find an {\it a priori} proof of the
conjecture failed.

\subsection{Standard filtrations}

Let $L$ be a simple Lie algebra over $F$ and $L_{(0)}$ a maximal
subalgebra of $L$. Let $L_{(-1)}$ be a subspace of $L$ such that
$L_{(0)}\subsetneq L_{(-1)}$ and $[L_{(0)},L_{(-1)}]\subset
L_{(-1)}$, and assume further that $L_{(-1)}/L_{(0)}$ is an
irreducible $L_{(0)}$-module. Following Weisfeiler \cite{We} we
define the {\it standard filtration} of $L$ associated with the pair
$(L_{(0)}, L_{(-1)})$ by setting
\begin{eqnarray*}
L_{(i+1)}&=&\{x\in L_{(i)}\,\vert\, [x,L_{(-1)}]\subset
L_{(i)}\},\qquad i\ge 0,\\
L_{(-i-1)}&=&[L_{(-i)},L_{(-1)}]+L_{(-i)},\qquad \qquad \quad i>0.
\end{eqnarray*}
Since $L_{(0)}$ is a maximal subalgebra of $L$ this filtration is
exhaustive. Since $L$ is simple, the filtration is separating. So
there are $s_1>0$ and $s_2\ge 0$ such that
$$L=L_{(-s_1)}\supset\ldots\supset L_{(0)}\supset\ldots\supset
L_{(s_2+1)}=(0).$$ By construction, all subspaces $L_{(i)}$ of $L$
are invariant under the action of the restricted subalgebra of
$\Der L$ generated by $\mbox{ad}\,L_{(0)}$. A standard filtration
is called {\it long} if $L_{(1)}\ne (0)$.

Now let $G\,=\,\bigoplus_{i\in\Z}\,G_i$ be a graded Lie algebra,
that is $[G_i,G_j]\subset G_{i+j}$ for all $i,j\in\Z$. The
following four conditions occur very frequently in Classification
Theory:
\begin{itemize}
\item[(g1)] $G_{-1}$ is an irreducible and faithful $G_0$-module;
\item[(g2)] $G_{-i}=[G_{-i+1},G_{-1}]$ for all $i\ge 1$;
\item[(g3)] if $x\in G_i$, $i> 0$, and $[x,G_{-1}]=(0)$, then
$x=0$;
\item[(g4)] if $x\in G_{-i}$, $i>0$, and $[x,G_k]=(0)$ for all $k>0$, then
$x=0$.
\end{itemize}
The graded Lie algebra
$\mbox{gr}\,L\,=\,\bigoplus_{i=-s_1}^{s_2}\mbox{gr}_i\, L,$ where
$\mbox{gr}_i\,L = L_{(i)}/L_{(i+1)}$, corresponding to the
standard filtration above satisfies the conditions (g1), (g2),
(g3). The quotient of $\mbox{gr}\,L$ by its largest ideal
contained in $\sum_{i<-1}\,\mbox{gr}_i\, L$ satisfies all four
conditions (g1) -- (g4).

\section{Classes of simple Lie algebras}

The main conjecture on the structure of finite dimensional simple
Lie algebras over algebraically closed fields of characteristic
$p$ is known as {\bf the generalized Kostrikin--Shafarevich
conjecture}. It states the following:

\smallskip

\noindent {\it For $p>5$, any finite dimensional simple Lie
algebra over $F$ is either classical or isomorphic to one of the
filtered Lie algebras of Cartan type.}

\smallskip

\noindent This conjecture is due to Kac \cite{K, K2} who
formulated it for $p>3$ (see also \cite{K71}). Our next goal is to
give a detailed description of the Lie algebras mentioned in the
generalized Kostrikin--Shafarevich conjecture.

\subsection{Classical Lie algebras}
Let $\mathfrak g$ be a simple Lie algebra over $\mathbb{C}$,
$\mathfrak h$ a Cartan subalgebra of $\mathfrak g$, $\Phi=\Phi
(\mathfrak{g}, \mathfrak{h})$ the corresponding root system, and
$\Delta=\{\alpha_1,\ldots,\alpha_l\}$ a basis of simple roots in
$\Phi$. For $\alpha,\beta\in \Phi$ set $\langle \beta,
\alpha^{\vee}\rangle=2(\beta\vert\alpha)/(\alpha\vert\alpha)$,
where, as usual, $(\,\cdot\,\vert\,\cdot\,)$ denotes a scalar
product on the $\mathbb{R}$-span of $\Phi$ invariant under the
Weyl group of $\Phi$.

\begin{theo}[\cite{Che}] The Lie algebra $\mathfrak g$ has a basis
$${\mathcal B}\,=\,\{e_{\alpha}\,\vert\,\alpha\in
\Phi\}\cup\{h_i\,\vert\,1\leq i\leq l\}$$ such that the following
conditions hold:
\begin{enumerate}
\item $[h_i,h_j]=0,\quad 1\leq i,j \leq l.$ \item $[h_i,
e_{\beta}]=\langle \beta, \alpha_{i}^{\vee}\rangle\,
e_{\beta},\quad 1\leq i \leq l,\quad \beta\in \Phi.$ \item
$[e_{\alpha},e_{-\alpha}]=h_{\alpha}$ is a $\Z$-linear combination
of $h_1,\ldots, h_l.$ \item Let $\alpha, \beta\in \Phi$, $\beta
\neq \pm\alpha$, and let $\{\beta-q\alpha,\ldots,\beta+r\alpha\}$
be the $\alpha$-string through $\beta$. Then
$[e_{\alpha},e_{\beta}]=0$ if $\alpha+\beta\not\in \Phi$ and
$[e_{\alpha},e_{\beta}]=\pm (q+1)e_{\alpha+\beta}$ if
$\alpha+\beta\in \Phi$. Moreover, $q\in\{0,1,2\}$ if
$\alpha+\beta\in\Phi$.
\end{enumerate}
\end{theo}
\noindent The $\Z$-span $\mathfrak{g}_{\Z}$ of $\mathcal B$ is a
$\Z$-form in $\mathfrak g$ closed under taking Lie brackets.
Therefore, $\mathfrak{g}_F:=\mathfrak{g}_{\mathbb
Z}\otimes_{\mathbb Z} F$ is a Lie algebra over $F$ with basis
${\mathcal B}\otimes 1$ and structure constants obtained from
those for $\mathfrak{g}_{\mathbb Z}$ by reducing modulo $p$. For
$p>3$, the Lie algebra $\mathfrak{g}_F$ fails to be simple if and
only if the root system $\Phi=\Phi(\mathfrak{g},\mathfrak{h})$ has
type $A_l$ where  $l=mp-1$ for some $m\in {\mathbb N}$. If $\Phi$
has type $A_{mp-1}$ then $\mathfrak{g}_F\cong \mathfrak{sl}(mp)$
has a one-dimensional center (consisting of scalar matrices) and
the Lie algebra $\mathfrak{g}_F/\mathfrak{z}(\mathfrak{g}_F)\cong
\mathfrak{psl}(mp)$ is simple. The simple Lie algebras over $F$
thus obtained are called {\em classical}.

All classical Lie algebras are restricted with $p$th power map
given by $(e_{\alpha}\otimes 1)^{[p]}=0$ and $(h_{i}\otimes
1)^{[p]}=h_i\otimes 1$ for all $\alpha\in \Phi$ and $1\le i\le l$.
As in characteristic $0$, they are parametrized by Dynkin diagrams
of types $A_n,\ B_n,\ C_n,\ D_n,\ G_2,\  F_4,\ E_6,\ E_7,\  E_8.$
We stress that, by abuse of characteristic $0$ notation, the
classical simple Lie algebras over $F$ include the Lie algebras of
simple algebraic $F$-groups of exceptional types. All classical
simple Lie algebras are closely related to simple algebraic groups
over $F$.

A Lie algebra $L$ of characteristic $p>3$ is called {\it almost
classical} if
$$\ad\mathfrak{g}\subset L \subset\mbox{Der}\,\mathfrak{g}$$
where $\mathfrak g$ is a direct sum of classical simple Lie
algebras. One of the examples of such algebras is the Lie algebra
$\mathfrak{pgl}(n):=\mathfrak{gl}(n)/F I_{n}$. When $p$ does not
divide $n$, we have that $\mathfrak{pgl}(n)\cong \mathfrak{sl}(n)$
as Lie algebras. However, $\mathfrak{pgl}(mp)\not\cong
\mathfrak{sl}(mp)$, because for $p>2$ the Lie algebra $
\mathfrak{sl}(mp)$ is perfect with a $1$-dimensional center, while
the Lie algebra $\mathfrak{pgl}(mp)$ is centerless and
$\big[\mathfrak{pgl}(mp),\mathfrak{pgl}(mp)\big]=\,
\mathfrak{psl}(mp)$ is an ideal of codimension $1$ in
$\mathfrak{pgl}(mp)$. It is easy to see that the Lie algebra
$\mathfrak{pgl}(mp)$ is almost classical.

All almost classical Lie algebras are semisimple, but the case of
$\mathfrak{pgl}(mp)$ shows that they are not always direct sums of
classical simple Lie algebras. Kostrikin conjectured in \cite{K62,
K71} that for $p>5$ a Lie algebra $L$ over $F$ is almost classical
if and only if ${\mathcal S}(L)=\{0\}$ (a closely related
conjecture can be found in the last section of \cite{KS66}).
Kostrikin's conjecture was proved in \cite{P86a} for $p>5$ and in
\cite{P86c} for $p=5$.

\subsection{Graded Lie algebras of Cartan type}

In \cite{KS69}, Kostrikin and Shafarevich gave a unified
description of a large class of nonclassical simple Lie algebras
over $F$. Their construction was motivated by classical work of E.
Cartan \cite{C09} on infinite dimensional, simple transitive
pseudogroups of transformations. To define finite dimensional
modular analogues of complex Cartan type Lie algebras Kostrikin
and Shafarevich replaced formal power series algebras over
$\mathbb C$ by divided power algebras over $F$.

Let $\mathbb{N}_{0}^m$ denote the additive monoid of all
$m$-tuples of nonnegative integers. For $\alpha,\beta\in
\mathbb{N}_{0}^m$ define
${\alpha\choose\beta}\,=\,{\alpha(1)\choose\beta(1)} \cdots
{\alpha(m)\choose\beta(m)}$ and
$\alpha!\,=\,\prod_{i=1}^m\alpha(i)!$. For $1\le i\le m$ set
$\epsilon_i\,=\,(\delta_{i1},\ldots,\delta_{im})$ and
$\underline{1}=\epsilon_1+\ldots+\epsilon_m$.

Give the polynomial algebra $F[X_1,\ldots,X_m]$ its standard
coalgebra structure (with all $X_i$ being primitive) and denote by
${\mathcal O}(m)$ the graded dual of $F[X_1,\ldots, X_m]$, a
commutative associative algebra over $F$. It is well-known (and
easily seen) that ${\mathcal O}(m)$  has basis
$\{x^\alpha\,\vert\,\alpha\in\mathbb{N}_0^m\}$ and the product in
${\mathcal O}(m)$ is given by $$x^\alpha
x^\beta={\alpha\choose\beta}x^{\alpha+\beta}\qquad \mbox{for all
}\ \ \alpha,\beta\in\mathbb{N}_0^m.$$ We write $x_i$ for
$x^{\epsilon_i}\in {\mathcal O}(m)$, $1\le i\le m$. For each
$m$-tuple $\underline{n}\in\mathbb{N}^m$ we denote by ${\mathcal
O}(m;\underline{n})$ the $F$-span of all $x^\alpha$ with $0\le
\alpha(i)<p^{n_i}$ for $i\le m$. This is a subalgebra of
${\mathcal O}(m)$ of dimension $p^{|\underline{n}|}$ where
$|\underline{n}|=n_1+\cdots+n_m$. Note that ${\mathcal
O}(m;\underline{1})$ is isomorphic to the truncated polynomial
algebra $F[X_1,\ldots, X_m]/(X_1^p,\ldots,X_m^p)$.

Assigning degree $|\alpha|=\alpha(1)+\cdots+\alpha(m)$ to each
$x^\alpha\in {\mathcal O}(m)$ gives rise to a grading of the
algebra ${\mathcal O}(m)$, called {\it standard}. Each ${\mathcal
O}(m;\underline{n})$ is a graded subalgebra of ${\mathcal O}(m)$.
The $k$th graded component of ${\mathcal O}(m)$ is denoted by
${\mathcal O}(m)_k$. The subspaces ${\mathcal
O}(m)_{(k)}\,:=\,\bigoplus_{i\geq k}\,{\mathcal O}(m)_i$ form a
decreasing filtration of ${\mathcal O}(m)$, called {\it the
standard filtration}. The completion of ${\mathcal O}(m)$ relative
to its standard filtration is denoted by ${\mathcal O}((m))$. The
elements of ${\mathcal O}((m))$ are the infinite formal sums of
the form $\sum_{\alpha}\lambda_{\alpha}\,x^{\alpha}$ with
$\lambda_{\alpha}\in F$. The algebra ${\mathcal O}((m))$ is
linearly compact and ${\mathcal O}(m)$ is canonically embedded
into ${\mathcal O}((m))$. The subspaces ${\mathcal
O}((m))_{(k)}:=\{\sum_{|\alpha|\geq
k}\lambda_{\alpha}\,x^\alpha\,\vert\,\lambda_\alpha\in F\}$ and
${\mathcal O}((m))_{k}:={\mathcal O}(m)_k$ induce a  decreasing
filtration and topological grading of ${\mathcal O}((m))$,
respectively. These are, again, called {\it standard}.

There is a family of continuous maps $\{y\mapsto y^{(s)}\,\vert\,
s\in\mathbb{N}_0\}$ from ${\mathcal O}((m))_{(1)}$ into ${\mathcal
O}((m))$, called {\it divided power maps}, such that
\begin{eqnarray*}
x^{(0)}\, =\, 1 &\qquad \ \ &\mbox{ for all   }\ \  x\in {\mathcal
O}((m))_{(1)};\\ (x^\alpha)^{(s)}\,=\, \Big( (s\alpha)!/(\alpha
!)^s s!\Big) x^{s\alpha} &\qquad \ \ &\mbox{ for all   }\ \ \alpha
\neq (0,\ldots, 0);\\ (\lambda x)^{(s)}\, =\,
\lambda^sx^{(s)}&\qquad \
\ & \mbox{ for all   }\ \ \lambda\in F, \ x\in {\mathcal O}((m))_{(1)};\\
(x+y)^{(s)}\,=\, \sum_{i=0}^{s} x^{(i)}y^{(s-i)}&\qquad \ \
&\mbox{ for all   }\ \ x,y\in {\mathcal O}((m))_{(1)}.
\end{eqnarray*}
A continuous automorphism $\phi$ (respectively, derivation $D$) of
the topological algebra ${\mathcal O}((m))$ is called {\it
admissible} (respectively, {\it special}) if $\phi(x^{(s)})= (\phi
x)^{(s)}$ (respectively, $D(x^{(s)})=x^{(s-1)} Dx$) for all $x\in
{\mathcal O}((m))_{(1)}$ and all $s\in\mathbb{N}_0$. For $1\le
i\le m$, the $i$th  partial derivative $\partial_i$ of ${\mathcal
O}((m))$ is defined as the special derivation of ${\mathcal
O}((m))$ with the property that
$\partial_{i}(x^\alpha)=x^{\alpha-\epsilon_i}$ if $\alpha(i)>0$
and $0$ otherwise. Each admissible automorphism of ${\mathcal
O}((m))$ respects the standard filtration of ${\mathcal O}((m))$.
Each finite dimensional subalgebra ${\mathcal O}(m;\underline{n})$
is stable under the partial derivatives
$\partial_1,\ldots,\partial_m$.

The set $W((m))$ of all special derivations of ${\mathcal O}((m))$
is an infinite dimensional Lie subalgebra of ${\rm Der}\,
{\mathcal O}((m))$ and an ${\mathcal O}((m))$-module, via
$(fD)(x)=f Dx$ for all $f\in {\mathcal O}((m))$ and $D\in W((m))$.
Since each $D\in W((m))$ is uniquely determined by its values
$Dx_1,\ldots, Dx_m$, the Lie algebra $W((m))$ is a free ${\mathcal
O}((m))$-module with basis $\partial_1,\ldots,\partial_m$. The
subspaces
$$W((m))_{k}\,:=\,\bigoplus_{i=1}^{m}{\mathcal O}((m))_{k+1}\partial_i \ \,
\mbox{ and }\ \ W((m))_{(k)}\,:=\,\bigoplus_{i=1}^{m}{\mathcal
O}((m))_{(k+1)}\partial_i
$$ for $k\ge -1$ form a topological grading and decreasing
filtration of $W((m))$, respectively. Needless to say, both are
called standard. Note that $$[W((m))_{(i)}, W((m))_{(j)}]\subset
W((m))_{(i+j)} \quad \ \mbox{for all }\ \ i\ge -1,\ j\ge 0.$$ The
group ${\rm Aut}_{c}\,{\mathcal O}((m))$ of all admissible
automorphisms acts on $W((m))$ by the rule $D\mapsto
D^\phi\,:=\,\phi^{-1}D\phi$, where $\phi\in{\rm Aut}_c\, {\mathcal
O}((m))$ and $D\in W((m))$, and respects the standard filtration
of $W((m))$.

The {\it general} Cartan type Lie algebra $W(m;\underline{n})$ is
the ${\mathcal O}(m;\underline{n})$-submodule of $W((m))$
generated by the partial derivatives
$\partial_1,\ldots,\partial_m$.  The Lie algebra
$W(m;\underline{n})$ is a subalgebra ${\rm Der}\, {\mathcal
O}(m;\underline{n})$. When $\underline{n}=\underline{1}$, it is
isomorphic to the full derivation algebra of $F[X_1,\ldots,
X_m]/(X_{1}^p,\ldots, X_{m}^p)$, a truncated polynomial ring in
$m$ variables. In the literature, $W(m;\underline{n})$ is often
referred to as a Lie algebra of {\it Witt} type.  Since
$W(m;\underline{n})$ is obviously a free ${\mathcal
O}(m;\underline{n})$-module of rank $m$, we have that $\dim\,
W(m;\underline{n})\,=\, mp^{|\underline{n}|}$. The Lie algebra
$W(m;\underline{n})$ is simple unless $(p,m)=(2,1)$. If
$\underline{n}\neq \underline{1}$ and $n_r\ne1$ then
$\partial_{r}^p\neq 0$ on ${\mathcal O}(m;\underline{n})$. Since
$\partial_r^p$ is not a special derivation of ${\mathcal O}((m))$
it follows that $W(m;\underline{n})$ is restrictable if and only
if $\underline{n}=\underline{1}$.

Give the ${\mathcal O}((m))$-module $$\Omega^1((m))\,:=\,{\rm
Hom}_{{\mathcal O}((m))}\Big(W((m)), {\mathcal O}((m))\Big)$$  a
$W((m))$-module structure by setting
$(D\alpha)(D')\,:=\,D(\alpha(D'))-\alpha([D,D'])$ for all $D,D'\in
W((m))$ and $\alpha\in\Omega^1((m))$, and define $d\,\colon\,
{\mathcal O}((m))\longrightarrow \Omega^1((m))$ by the rule
$(df)(D)=Df$ for all $D\in W((m))$ and $f\in {\mathcal O}((m))$.
Notice that $d$ is a homomorphism of $W((m))$-modules and
$\Omega^1((m))$ is a free ${\mathcal O}((m))$-module with basis
$dx_1,\ldots, dx_m$. Let
$$\Omega((m))\,=\,\bigoplus_{0\le k\le m}\,\Omega^k((m))$$ be the
exterior algebra, over ${\mathcal O}((m))$, on $\Omega^1((m))$.
Then $\Omega^0((m))={\mathcal O}((m))$ and each graded component
$\Omega^k((m))$, $k\ge 1$, is a free ${\mathcal O}((m))$-module
with basis $\{dx_{i_1}\wedge\ldots\wedge dx_{i_k}\,\vert\,1\le
i_1<\ldots< i_k\le m\}$. The elements of $\Omega((m))$ are called
{\it differential forms} on ${\mathcal O}((m))$.

The map $d$ extends (uniquely) to a zero-square linear operator of
degree $1$ on $\Omega((m))$ such that $$
d(f\omega)\,=\,(df)\wedge\omega+fd(\omega),\ \ \
d(\omega_1\wedge\omega_2)\,=\,d(\omega_1)\wedge\omega_2+
(-1)^{\deg(\omega_1)}\omega_1\wedge d(\omega_2)$$ for all $f\in
{\mathcal O}((m))$ and all homogeneous
$\omega,\omega_1,\omega_2\in\Omega((m))$. For $D\in W((m))$, we
have that $D(f\omega)=(Df)\omega+fD(\omega)$. It follows that each
$D\in W((m))$ extends to a derivation of the $F$-algebra
$\Omega((m))$. All such derivations commute with $d$. The group
${\rm Aut}_{c}\,{\mathcal O}((m))$ acts on $\Omega^1((m))$ by the
rule
$$(\phi \omega)(D)\,:=\,\phi(\omega(D^\phi))$$ for all
$\phi\in{\rm Aut}_c\,{\mathcal O}((m))$, $\omega\in\Omega^1((m))$,
$D\in W((m))$. Moreover, $$\phi(f\omega)\,=\,\phi(f)\phi(\omega) \
\mbox{ and } \ \phi\circ d\,=\,d\circ\phi$$ for all $\phi\in{\rm
Aut}_c\,{\mathcal O}((m))$, $\omega\in\Omega((m))$, $f\in
{\mathcal O}((m))$. It follows that the action of ${\rm
Aut}_c\,{\mathcal O}((m))$ on $\Omega^1((m))$ extends to an
embedding ${\rm Aut}_c\,{\mathcal O}((m))\hookrightarrow {\rm
Aut}_{F}\,\Omega((m)).$ It can be shown that
$$D^\phi(\omega)\,=\,\phi^{-1}(D(\phi(\omega))$$ for all $D\in
W((m))$, $\phi\in{\rm Aut}_c\,{\mathcal O}((m))$,
$\omega\in\Omega((m))$.

Each $m$-tuple $\underline{r}$ of nonnegative integers induces a
grading of the algebra ${\mathcal O}(m)$ defined by assigning
${\rm deg}\, (x^\alpha)=r(1)\alpha (1)+\cdots+r(m)\alpha(m)$ to
each monomial $x^\alpha\in {\mathcal O}(m)$. Such a grading, in
turn, induces (topological) gradings and decreasing filtrations of
the algebras ${\mathcal O}(m;\underline{n})$, ${\mathcal O}((m))$,
$W(m;\underline{n})$, and $W((m))$. It also induces a topological
grading of the algebra $\Omega((m))$ which extends that of
${\mathcal O}((m))=\Omega^0((m))$. The differential $d$ of
$\Omega((m))$ preserves all components of this grading. The
gradings and filtrations thus obtained are all said to be of {\it
type} $\underline{r}$. In this new terminology, the standard
gradings and filtrations defined above are all of type
$\underline{1}$.

As in the characteristic $0$ case, the three differential forms
below are of particular interest: $$\begin{array}{rcll} \omega_S &
:= & dx_1\wedge\ldots\wedge dx_m, & m\geq 3,\\ \omega_H & := &
\sum_{i=1}^r dx_i\wedge dx_{i+r}, & m=2r\ge 2,\\ \omega_K & := &
dx_{2r+1}+\sum_{i=1}^r(x_{i+r}dx_i-x_idx_{i+r}) ,& m=2r+1\ge 3.
\end{array}
$$
These forms give rise to the following Lie algebras:
\begin{eqnarray*}
S((m)) & := & \{D\in W((m))\mid D(\omega_S)=0\},\\ &&\hspace{2cm}
special\ Lie\ algebra,\\[0.3ex] H((m)) & := & \{D\in W((m))\mid
D(\omega_H)=0\},\\ &&\hspace{2cm} Hamiltonian\ Lie\
algebra,\\[0.3ex] K((m))& := & \{D\in W((m))\mid D(\omega_K)\in
{\mathcal O}((m))\omega_K\},\\ &&\hspace{2cm} contact\ Lie\
algebra.
\end{eqnarray*}
\noindent
Define Lie algebras $CS((m))$ and $CH((m))$ by setting
\begin{eqnarray*}
CS((m))&:=&\{D\in W((m)\mid D(\omega_S)\in F\omega_S\},
\\ CH((m))&:=&\{D\in W((m))\mid D(\omega_H)\in F\omega_H\}.
\end{eqnarray*}
\noindent Obviously, $CX((m))^{(1)}\subseteq X((m))$ for
$X\in\{S,H\}$. For $X\in\{W, S, CS, H, CH\}$, set
$\underline{r}_X=\epsilon_1+\cdots+\epsilon_m={\underline 1}$. For
$X=K$, set $\underline{r}_X=\epsilon_1+\cdots+
\epsilon_{m-1}+2\epsilon_m={\underline 1}+\epsilon_m$. For
$X\in\{W,S,CS,H,CH,K\}$ and ${\underline n}\in\mathbb{N}^m$,
define
$$X(m;{\underline n})\,=\,X((m))\cap W(m;{\underline n}).$$ Each
$X(m;{\underline n})$ is a graded subalgebra of the Lie algebra
$X((m))$ regarded with its grading of type ${\underline r}_X$. The
graded components of $X(m;{\underline n})$ are denoted by
$X(m;{\underline n})_i\,$, $\,i\in\Z$. Note that $X(m;{\underline
n})_i\,=\,(0)$ for $i\le -2$ if $X\neq K$. Also,
$\dim\,K(m;{\underline n})_{-2}=1$ and $K(m;{\underline
n})_{i}\,=\, (0)$ for $i\le -3$.

Suppose $p>3$. In \cite{KS69}, it was shown that the Lie algebras
$S(m;{\underline n})^{(1)}$, $H(m;{\underline n})^{(1)}$ and
$K(m;{\underline n})^{(1)}$ are simple for $m\ge 3$ and that so is
$H(2;{\underline n})^{(2)}$. Moreover, $K(m;{\underline n})\,=\,
K(m;{\underline n})^{(1)}$ unless $p\,\vert\,(m+3)$. For
$X\in\{W,S,CS,H,CH,K\}$ any $\underline{r}_X$-graded Lie
subalgebra of $X(m;{\underline n})$ containing $X(m;{\underline
n})^{(\infty)}$ is called a finite dimensional {\it graded} Lie
algebra of Cartan type. According to \cite{KS66} the Lie algebra
$X(m;{\underline n})^{(\infty)}$ is restrictable if and only if
${\underline n}={\underline 1}$.

{\bf The original Kostrikin--Shafarevich conjecture} \cite{KS66}
of 1966 states the following:

\smallskip

\noindent {\it For $p>5$, any finite dimensional restrictable
simple Lie algebra over $F$ is either classical or isomorphic to
one of the Lie algebras $W(m;{\underline 1}),$ $m\ge 1\,$,
$S(m;{\underline 1})^{(1)},$ $m\ge 3\,$, $H(m;{\underline
1})^{(2)},$ $m\ge 2$, $\,K(m;{\underline 1})^{(1)},$ $m\ge 3$}.

\subsection{Filtered Lie algebras of Cartan type}

In order to give a unified description of all known finite
dimensional simple Lie algebras of characteristic $p>5$ Kac
\cite{K2} and Wilson \cite{Wil69, Wil76} introduced certain
filtered deformations of finite dimensional graded Lie algebras of
Cartan type. A streamlined treatment of these algebras is given
\cite{St04}.

We first outline Wilson's original approach. Let
$X\in\{W,S,H,K\}$, $\underline{n}\in\mathbb{N}_0^m$, and let
$\Phi$ be an admissible automorphism of ${\mathcal O}((m))$. For
$X=K$ assume further that $\Phi$ respects the
$\underline{r}_X$-filtration of ${\mathcal O}((m))$ (if $X\neq K$
this assumption is fulfilled automatically). Define
$$X(m;\underline{n};\Phi)\ :=\ (\Phi^{-1}\circ X((m))\circ
\Phi)\cap W(m;\underline{n}).$$ It is clear from the definition
that $X(m;\underline{n}; \mbox{Id})\,=\, X(m;\underline{n})$ and
$W(m;\underline{n};\Phi)\,=\,W(m;\underline{n})$.

\begin{defi}[\cite{Wil76}] The Lie algebra $X(m;\underline{n};\Phi)^{(\infty)}$ is called a
{\rm filtered Lie algebra of Cartan type} if $X(m;n;\Phi)$
satisfies the following two conditions:
\begin{itemize}
\item[1.] $X(m;\underline{n};\Phi)\cap
W(m;\underline{n})_{(2+\delta_{X,K}),X}\,\neq \,(0);$ \item[2.]
$X(m;\underline{n};\Phi)+(\Phi\circ X((m))\circ\Phi^{-1})\cap
W(m;\underline{n})_{(1+\delta_{X,K}),X}\,=\,\Phi\circ
X((m))\circ\Phi^{-1}.$
\end{itemize}
Here $W(m;\underline{n})_{(k),X}$ denotes the $k${\rm th}
component of the $\underline{r}_X$-filtration of
$W(m;\underline{n})$.
\end{defi}
\noindent The embedding of a filtered Cartan type Lie algebra
$X(m;\underline{n};\Phi)^{(\infty)}$ into the Lie algebra
$W(m;\underline{n})$ regarded with its filtration of type $r_X$,
induces a natural filtration of
$X(m;\underline{n};\Phi)^{(\infty)}$. The corresponding graded
algebra is isomorphic to a graded Cartan type Lie algebra
(possibly of type $CS$ or $CH$) containing
$X(m;\underline{n})^{(\infty)}$ as a minimal ideal. The subalgebra
$L_{(0)}\,=\,X(m;\underline{n};\Phi)^{(\infty)} \cap
W(m;\underline{n})_{(0)}$ is called the {\it standard maximal
subalgebra} of $L\,= \,X(m;\underline{n};\Phi)^{(\infty)}$. For
$p>3$, this subalgebra can be characterized as the unique proper
subalgebra of maximal dimension in $L$; see \cite{Kr}, \cite{Skk},
\cite{St04}. As a consequence, $L_{(0)}$ is stable under all
automorphisms of $L$. For $p>3$, each Cartan type Lie algebra $L$
is simple (\cite{Wil76}).

The following important abstract characterization of filtered
Cartan type Lie algebras is due to Wilson \cite{Wil76}. Let
$\mathcal L$ be a Lie algebra over $F$ and let ${\mathcal L}_0$ be
a subalgebra of $\mathcal L$. Then we have a natural
representation $\rho\,\colon\, {\mathcal L}_0\rightarrow
\mathfrak{gl}({\mathcal L}/{\mathcal L}_0)$ of the Lie algebra
${\mathcal L}_0$ given by
$$(\rho(x))(y+{\mathcal L}_0)=[x,y]+{\mathcal L}_0\ \ \ \mbox{for all }\,
x\in {\mathcal L}_0, \ y\in{\mathcal L}.$$
\begin{theo}[Wilson's Theorem] Let $\mathcal L$ be a
simple Lie algebra over $F$ and suppose that ${\mathrm char}\,
F=p>3$. Then $\mathcal L$ is isomorphic to a finite dimensional
filtered Cartan type Lie algebra if and only if $\mathcal L$ is
strongly degenerate and contains a maximal subalgebra ${\mathcal
L}_0$ such that either ${\mathcal L}_0$ has codimension $1$ in
$\mathcal L$ or else $\rho({\mathcal L}_0)$ contains a linear
transformation $Y$ of rank $1$ such that $[Y,[Y,\rho({\mathcal
L}_0)]]\neq (0)$.
\end{theo}

Kac's approach \cite{K2} to filtered Cartan type Lie algebras
pushed further by Skryabin in \cite{Sk', Sk'', Sk, Skk', Skk}
involved more general differential forms in $\Omega((m))$.
Combined with Wilson's theorem it eventually led to a complete
classification of filtered Lie algebras of Cartan type.

Recall that the algebra ${\mathcal O}((m))$ is linearly compact.
Given a unital associative subalgebra $B$ of ${\mathcal O}((m))$
we let $W(B)$ and $W(B)_{(0)}$ denote the normalizers of $B$ and
$B\cap {\mathcal O}((m))_{(0)}$ in $W((m))$, respectively. We
denote by $B^*$ the group of invertible elements of $B$. Following
\cite{Sk} we say that $B$ is an {\it admissible} subalgebra of
${\mathcal O}((m))$ if $B$ is closed in ${\mathcal O}((m))$ and
$W(B)_{(0)}$ has codimension $m$ in $W(B)$. This definition is
inspired by a crucial definition in \cite{K2}. Any finite
dimensional subalgebra of ${\mathcal O}((m))$ of the form $\phi
({\mathcal O}(m;\underline{n}))$ with $\phi\in\mbox{Aut}_{c}\,
{\mathcal O}((m))$ and $\underline{n}\in\mathbb{N}^m$ is
admissible. Conversely, given a finite dimensional admissible
subalgebra $B\subset {\mathcal O}((m))$ there are an automorphism
$\phi\in\mbox{Aut}_{c}\, {\mathcal O}((m))$ and a tuple
$\underline{n}\in \mathbb{N}^m$ such that $B\,=\,\phi({\mathcal
O}(m;\underline{n}))$; see \cite{Sk}.

To ease notation we set $\Omega=\Omega((m))$,
$\Omega^k=\Omega^k((m))$ for $0\le k\le m$, and put $\Omega^{\,\rm
even}:=\bigoplus_{i\ge 0}\,\Omega^{2i}$. Observe that
$\Omega^{\,\rm even}$ is a commutative algebra over $F$ and
$W((m))$ acts on $\Omega^{\,\rm even}$ as derivations. The
subspace $\Omega^{\,\rm even}_+:={\mathcal
O}((m))_{(1)}\oplus\bigoplus_{i\ge 1}\,\Omega^{2i}$ is a maximal
ideal of $\Omega^{\,\rm even}$ which intersects trivially with
$(\Omega^{\,\rm even})^{W((m))}=F1$. It is well-known that the
first cohomology group $H^1\big(W((m)),\,\Omega^{\,\rm even}\big)$
vanishes; see \cite[Theorem~7.5]{Sk} for example. According to
\cite[Proposition~1.2]{Sk} this implies that there exists unique
system of divided powers $\omega\mapsto \omega^{(s)},\,s\ge 0,$ on
$\Omega_+^{\,\rm even}$ with respect to which $W((m))$ acts on
$\Omega^{\,\rm even}$ as special derivations. It has the property
that $\omega^{(s)}\in\Omega^{2is}$ whenever $\omega\in\Omega^{2i}$
and $s\ge 1$.

Recall that a differential form $\omega\in\Omega$ is called {\it
closed} if $d\omega=0$. Following \cite{K2} we say that
$\omega\in\Omega^m$ is {\it nondegenerate} if $m\ge 2$ and
$\omega=\varphi\, dx_1\wedge\ldots\wedge dx_m$ for some
$\varphi\in {\mathcal O}((m))^*$. We call $\omega\in\Omega^2$ {\it
nondegenerate} if $m=2r\ge 2$, $\omega$ is closed, and the form
$\omega^{(r)}\in\Omega^m$ is nondegenerate (if $m=2$ this is
consistent with the previous definition). Finally, we say that
$\omega\in\Omega^1$ is {\it nondegenerate} if $m=2r+1\ge 3$ and
$(d\omega)^{(r)}\wedge\omega\in\Omega^m$ is nondegenerate.

Given a finite dimensional admissible subalgebra $B$ of ${\mathcal
O}((m))$ we let $\Omega(B)\,=\,\bigoplus_{k=0}^m\,\Omega^k(B)$
denote the $B$-subalgebra of $\Omega$ generated (over $B$) by
$dB$. For $f\in {\mathcal O}((m))_{(1)}$ we set  $\exp f:=
\sum_{i\ge 0}\, f^{(i)},$ an element in ${\mathcal O}((m))^*$. Let
${\bf s}(B)$ (respectively, ${\bf h}(B))$ denote the set of all
nondegenerate forms $\omega\in \Omega^m$ (respectively, $\omega\in
\Omega^2$) such that $\omega=(\exp u)\,\omega'$ for some
$\omega'\in\Omega(B)$ and $u\in {\mathcal O}((m))_{(1)}$
satisfying $du\in\Omega^1(B)$. Let ${\bf k}(B)$ denote the set of
all nondegenerate forms in $\Omega^1(B)$.

For $\omega\in {\bf s}(B)$ define the Lie algebras
\begin{eqnarray*}
S(B,\omega)&:=&\{D\in W(B)\,\vert\,
D\omega=0\};\\
CS(B;\omega)&:=&\{D\in W(B)\,\vert\,D\omega\in F\omega\}.
\end{eqnarray*}
For $\omega\in{\bf h}(B)$ define the Lie algebras
\begin{eqnarray*}
H(B;\omega)&:=&\{D\in
W(B)\,\vert\,D\omega=0\};\\
CH(B;\omega)&:=&\{D\in W(B)\,\vert\,D\omega\in F\omega\}.
\end{eqnarray*}
For $\omega\in{\bf k}(B)$ define the Lie algebra
\begin{eqnarray*}
K(B;\omega)&:=&\{D\in
W(B)\,\vert\,D\omega\in
B\,\omega\}.
\end{eqnarray*}
It is proved in \cite{K2, Skk', Skk} that except for two cases in
characteristic $2$ the Lie algebras $$W(B),\, S(B;\omega)^{(1)},\,
H(B;\omega)^{(2)},\, K(B;\omega)^{(1)}$$ are simple. Dimensions
and explicit bases of the Lie algebras $X(B;\omega)$,
$X(B;\omega)^{(1)}$ and $X(B;\omega)^{(2)}$ are found in
\cite{Skk, Kir1, Kir2} (see also \cite{BKK}). The most accessible
reference, by far, is \cite[Sect.~6]{St04}.

Any simple Lie algebra $L=X(B;\omega)^{(\infty)}$ is naturally
filtered and $\mbox{gr}\, L$, the corresponding graded algebra, is
isomorphic to a graded Lie algebra of Cartan type. In view of
Wilson's theorem this implies that for $p>3$ each
$X(B;\omega)^{(\infty)}$ is isomorphic to a filtered Cartan type
Lie algebra. The converse is also true: for $p>3$ any filtered
Cartan type Lie algebra $X(m;\underline{n};\Phi)^{(\infty)}$ is
isomorphic to one of $X(B;\omega)^{(\infty)}$ where
$B=\phi({\mathcal O}(m;\underline{n}))$ for some $\phi\in
\mbox{Aut}_{c}\,{\mathcal O}((m))$ (see \cite{K2}, \cite{Ku},
\cite{Sk}).

The $p$-structure of filtered Cartan type Lie algebras is described by
the following theorem.
\begin{theo} [\cite{K2, Skk}] Let $B=\phi({\mathcal O}(m;\underline{n}))$ where
$\phi\in{\rm Aut}_c\,{\mathcal O}((m))$.
\begin{enumerate}
\item
The Lie algebras $W(B)$, $CS(B;\omega)$, $CH(B;\omega)$,
$K(B;\omega)$ and $K(B;\omega)^{(1)}$ are restrictable if and only
if $\underline{n}=\underline{1}$.
\item
The Lie algebras $S(B;\omega)$ and $H(B;\omega)$ are
restrictable if and only if $\underline{n}=\underline{1}$
and $\omega\in\Omega(B)$.
\item
The Lie algebras $S(B;\omega)^{(1)}$, $H(B;\omega)^{(1)}$ and $H(B;\omega)^{(2)}$
are restrictable if and
only if $\underline{n}=\underline{1}$ and $\omega\in d\,\Omega(B)$.
\end{enumerate}
\end{theo}
\noindent It follows from Theorem 5 that for $p>3$ the Lie algebra
$X(m;\underline{n};\Phi)^{(\infty)}$ is restrictable if and only
if it is isomorphic to one of $W(m;\underline{1})$,
$S(m;\underline{1})^{(1)}$, $H(m;\underline{1})^{(2)}$,
$K(m;\underline{1})^{(1)}$.

The realizations of filtered Cartan type Lie algebras just
described are very useful in view of Kac's Isomorphism Theorem
which was later refined by Skryabin; see \cite{K2, Sk, Skk}. Let
$B$ (respectively, $B'$) be an admissible subalgebra of ${\mathcal
O}((m))$ (respectively, ${\mathcal O}((m'))$), and
$X,X'\in\{W,S,H,K\}$. Slightly abusing notation we set
$W(B;\omega)\,=\,W(B)$ and likewise for $W(B')$. We call a linear
map $$\sigma\,\colon\, X(B;\omega)^{(\infty)}\longrightarrow
X'(B';\omega')^{(\infty)}$$ {\it standard} if $\sigma(D)=\psi\circ
D\circ\psi^{-1}$ for all $D\in X(B;\omega)^{(\infty)}$, where
$\psi\,\colon\,{\mathcal O}((m))\stackrel{\sim}{\longrightarrow}
{\mathcal O}((m'))$ is a continuous isomorphism of divided power
algebras satisfying $\psi(B)=B'$ and $\psi(\omega)=C\omega'$ with
$C\in F^*$ for $\omega\in {\bf s}(B)\cup {\bf h}(B)$ and $C\in
{B'}^*$ for $\omega \in {\bf k}(B)$. Clearly, any standard map is
a Lie algebra isomorphism. Also, if $\sigma\,\colon\,
X(B;\omega)^{(\infty)}\rightarrow X'(B';\omega')^{(\infty)}$ is a
standard map then necessarily $m=m'$ and $X=X'$.
\begin{theo}[Isomorphism Theorem] Let $B,B'$ and
$X, X'$ be as above. Then with eight exceptions in characteristic $2$ and
three exceptions in characteristic
$3$ any isomorphism between the Lie algebras $X(B;\omega)^{(\infty)}$
and $X'(B';\omega')^{(\infty)}$ is standard.
\end{theo}
\noindent In our further discussion of the Lie algebras
$X(B;\omega)$ we shall assume (without loss of generality) that
$B={\mathcal O}(m;\underline{n})$. In this special case,
$X(B;\omega)$ is denoted by $X(m;\underline{n};\omega)$. The
corresponding set ${\bf x}(B)$ of nondegenerate forms will be
denoted by ${\bf x}(m;\underline{n})$. We shall also assume (as we
may) that $\underline{n}$ is a partition of $|\underline{n}|$,
that is $n_1\ge\ldots\ge n_m$. Let $G(m;\underline{n})$ denote the
set-wise stabilizer of ${\mathcal O}(m;\underline{n})$ in
$\mbox{Aut}_{c}\, {\mathcal O}((m))$. This is a connected
algebraic group with a large unipotent radical; see \cite{Wil71}
for more detail. It follows from Theorem 5 that with a few
exceptions in characteristics $2$ and $3$ the Lie algebras
$X(m;\underline{n};\omega)^{(\infty)}$ and
$X'(m';\underline{n}';\omega')^{(\infty)}$ are isomorphic if and
only if $m=m'$, $\underline{n}=\underline{n}'$ and
$g\omega=C\omega'$ for some $g\in G(m;\underline{n})$, where $C\in
F^*$ for $X\in\{S,H\}$ and $C\in {\mathcal O}(m;\underline{n})^*$
for $X=K$. Let $\mathfrak{S}({\underline{n}})$ denote the group of
all permutations $\pi$ of $\{1,2,\ldots,m\}$ such that $n_{\pi
i}=n_{i}$ for all $i$.

The orbits of $G(m;\underline{n})$ on ${\bf s}(m;\underline{n})$
are described in \cite{Tyu} and \cite{Wil80}. Let
$I(\underline{n})$ denote the subset of $\{1,2,\dots,m\}$
consisting of $1$ and all $k$ with $n_k<n_{k-1}$. Set
$\delta_{\underline{n}}=(p^{n_1}-1,\ldots, \,p^{n_m}-1)$.
According to \cite{Tyu, Wil80}, each $\omega\in{\bf
s}(m;\underline{n})$ is conjugate under $G(m;\underline{n})$ to a
nonzero scalar multiple of precisely one form in the set $$\{(\exp
x_i)\,\omega_S\,\vert\,i\in I(\underline{n})\}\cup\{\omega_S,\,
(1-x^{\delta_{\underline{n}}})\,\omega_S\}.$$ As a consequence,
for $p>2$ and $\underline{n}\in\mathbb{N}^m$ fixed, there are only
finitely many filtered Lie algebras of type
$S(m;\underline{n};\Phi)^{(\infty)}$ up to isomorphism.

The orbits of $G(m;\underline{n})$ on ${\bf k}(m;\underline{n})$
are described in \cite{KK$'$} in the simplest case
$\underline{n}=\underline{1}$ and in \cite{Sk'} for any
$\underline{n}$. Let ${\mathcal D}_{\bf k}$ denote the set of all
decompositions of $\{1,2,\ldots, m\}$ into a disjoint union of the
form $${\bf I}\,=\,\{i_0\}\sqcup\{i_1,i'_1\}\sqcup\ldots\sqcup
\{i_r,i'_r\},\ \ \ \ \ \ i_k<i'_k, $$ (different orderings of the
subsets within the union are not distinguished). The group
$\mathfrak{S}(\underline{n})$ acts on the set ${\mathcal D}_{\bf
k}$. Given ${\bf I}\in{\mathcal D_{\bf k}}$ define
$$\omega_{K,{\bf
I}}\,:=\,dx_{i_0}+\sum_{k=1}^{r}x_{i_k}dx_{i'_k},$$ an element in
${\bf k}(m;\underline{n})$. It is proved in \cite{Sk'} that for
$p>2$ each $\omega\in{\bf k}(m;\underline{n})$ is conjugate under
$G(m;\underline{n})$ to $f\omega_{K,{\bf I}}$ for some $f\in
{\mathcal O}(m;\underline{n})^*$ and ${\bf I}\in{\mathcal D}_{\bf
k}$. Moreover, the orbit of $f\omega_{K,{\bf I}}$ under
$G(m;\underline{n})$ intersects with ${\mathcal
O}(m;\underline{n})^*\omega_{K,{\bf I}'}$ for ${\bf
I}'\in{\mathcal D}_{\bf k}$ if and only if there is a
$\pi\in\mathfrak{S}(\underline{n})$ such that $\pi({\bf I})={\bf
I}'$.

Thus for $p>2$ any filtered Cartan type Lie algebra
$K(m;\underline{n};\omega)^{(\infty)}$ is isomorphic to a graded
Cartan type Lie algebra $K(m;\underline{n}')^{(1)}$ (here
$|\underline{n}'|=|\underline{n}|$ but in general $\underline{n}'$
need not be a partition of $|\underline{n}|$). It follows that for
$p>2$ and $\underline{n}\in\mathbb{N}^m$ fixed, there are only
finitely many Cartan type Lie algebras
$K(m;\underline{n};\Phi)^{(\infty)}$ up to isomorphism.

The orbit set ${\bf h}(m;\underline{n})/G(m;\underline{n})$ is
studied in \cite{K2, KK, BGOSW, Sk', Sk''}. In the simplest case
$\underline{n}=\underline{1}$ it is described in \cite{KK}. For an
arbitrary $\underline{n}$, a reasonably small set of
representatives for each $G(m;\underline{n})$-orbit in ${\bf
h}(m;\underline{n})$ is found in \cite{BGOSW}. A complete
description of ${\bf h}(m;\underline{n})/G(m;\underline{n})$ is
given in \cite{Sk', Sk''}.

Set ${\bf h}_1(m;\underline{n})\,:=\, {\bf
h}(m;\underline{n})\cap\Omega({\mathcal O}(m;\underline{n}))$ and
${\bf h}_2(m;\underline{n})\,:=\,{\bf h}(m;\underline{n})\setminus
{\bf h}_1(m;\underline{n}).$ Both ${\bf h}_1(m,\underline{n})$ and
${\bf h}_2(m;\underline{n})$ are $G(m;\underline{n})$-stable. The
orbit sets ${\bf h}_{2}(m;\underline{n})/G(m;\underline{n})$ and
${\bf k}(m+1;\underline{n})/\big(G(m+1;\underline{n})\ltimes
{\mathcal O}(m+1;\underline{n})^*\big)$ are somewhat similar to
each other. Let ${\mathcal D}_{\bf h}$ denote the set of all
decompositions
$${\bf I}\,=\,\{i_1,i'_1\}\sqcup\ldots\sqcup\{i_r,i'_r\},\ \ \ \ \
\ i_k<i'_k,$$ of $\{1,2,\ldots,m\}$ into a disjoint union of pairs
(different orderings of the pairs within the union are not
distinguished). The group $\mathfrak{S}(\underline{n})$ acts on
the set ${\mathcal D}_{\bf h}$. Given $i\in\{1,2,\ldots, m\}$ and
${\bf I}\in{\mathcal D_{\bf h}}$ define $$\omega_{H,i,{\bf
I}}\,:=\,d\Big(\exp
x_{i}\,\sum_{k=1}^{r}\,x_{i_k}dx_{i'_k}\Big),$$ an element in
${\bf h}_2(m;\underline{n})$. It is proved in \cite{Sk', Sk''}
that for $p>2$ each $\omega\in{\bf h}_2(m;\underline{n})$ is
conjugate under $G(m;\underline{n})$ to $\omega_{H,i,{\bf I}}$ for
some ${\bf I}\in{\mathcal D}_{\bf h}$ and $i\in\{1,2,\ldots, m\}$.
Moreover, the $G(m;\underline{n})$-orbit of $\omega_{H,i,{\bf I}}$
intersects with $F^*\omega_{H,j,{\bf I}'}$ for ${\bf
I}'\in{\mathcal D}_{\bf h}$ and $j\in\{1,2,\ldots, m\}$ if and
only if there is a $\pi\in\mathfrak{S}(\underline{n})$ such that
$\pi({\bf I})={\bf I}'$ and $\pi i=j$. As a consequence,  for
$p>2$ and $\underline{n}\in \mathbb{N}^m$ fixed, there are only
finitely many isomorphism classes of Lie algebras of the form
$H(m;\underline{n};\omega)^{(\infty)}$ with $\omega\in{\bf
h}_2(m;\underline{n})$.

The orbit set ${\bf h}_{1}(m;\underline{n})/G(m;\underline{n})$ is
{\it much} more complicated. It is no longer discrete, for $m\ge
4$, and this allows one to exhibit multiparameter families of
pairwise nonisomorphic simple Lie algebras of dimensions
$p^{|\underline{n}|}-2$ and $p^{|\underline{n}|}-1$. This
phenomenon was first discovered by Kac who disproved an earlier
conjecture of Kostrikin stating that no such families could exist
for $p>3$ (see \cite{K71}).

Let $J_{l}(\lambda)$ denote the Jordan block of order $l$ with
eigenvalue $\lambda\in F$. Let $O_l$ (respectively, $E_l$) denote
the zero (respectively, identity) matrix of order $l$. Let
\[C_l\,=\,\left[
\begin{array}{cccc}
0&1&\cdots&0\\
\vdots&\vdots&\ddots&\vdots\\
0&0&\cdots&1\\
1&0&\cdots&0
\end{array}
\right],
\]
a monomial matrix of order $l$. Given $d,s\in\mathbb{N}$ and
$\lambda\in F$ define $$C_{d,s}(\lambda)\,=\,\left[
\begin{array}{cccc}
O_s&E_s&\cdots&O_s\\
\vdots&\vdots&\ddots&\vdots\\
O_s&O_s&\cdots&E_s\\
J_{s}(\lambda)&O_s&\cdots&O_s
\end{array}
\right], $$ a block-monomial matrix of order $ds$. Let ${\mathcal
H}_m$ denote the set of all pairs of block-diagonal, skew-symmetric
matrices
$$(A,B)=\Big(\mbox{diag}(A_1,\ldots,A_k),\,\mbox{diag}
(B_1,\ldots,B_k)\Big)$$ of order $m=2r=2r_1+\ldots+2r_k\ge 2$ such
that $$A_{i}\,=\,\left[\begin{array}{cc}O_{r_i}&E_{r_i}\\
-E_{r_i}&O_{r_i}
\end{array}\right]$$ and
$B_i$ is one of
$$\left[\begin{array}{cc}O_{r_i}&J_{r_i}(0)\\-J_{r_i}(0)&O_{r_i}\end{array}\right],\
\, \left[\begin{array}{cc}O_{r_i}&
C_{d_i,s_i}(\lambda)\\-C_{d_i,s_i}(\lambda)&O_{r_i}\end{array}\right],
\ \,
\left[\begin{array}{cc}O_{r_i}&C_{r_i}\\-C_{r_i}&O_{r_i}\end{array}\right],$$
where $\lambda\ne 0$ and $r_i=d_is_i$. To each
$(A,B)=\big((a_{ij}),\,(b_{ij})\big)\in {\mathcal H}_m$ one
associates a differential form $\omega_{A,B}\in\Omega((m))$ by
setting
$$\omega_{A,B}\,=\,
\sum_{i<j}\,\Big(a_{ij}+b_{ij}x_{i}^{(p^{n_i}-1)}x_{j}^{(p^{n_j}-1)}\Big)dx_{i}\wedge
dx_j.$$ It is straightforward to see that $\omega_{A,B}\in{\bf
h}_1(m;\underline{n})$. One of the main results in \cite{Sk'} (see
also \cite{Sk''}) says that any $\omega\in{\bf
h}_1(m;\underline{n})$ is conjugate under $G(m;\underline{n})$ to
one of $\omega_{A,B}$ with $(A,B)\in{\mathcal H}_m$ (this holds in
all prime characteristics). Skryabin also found a necessary and
sufficient condition for two forms $\omega_{A,B}$ and
$\omega_{A',B'}$ to be conjugate under $G(m;\underline{n})$. It
involves an equivalence relation on the set of all pairs of
sequences of natural numbers, finite of equal length or periodic;
see \cite{Sk''} for more detail.

\subsection{Melikian algebras and their
relatives}

In this subsection we assume that  $p\in\{2,3,5\}$. Around 1980,
Melikian (a PhD student of Kostrikin at the time) discovered a new
series of finite dimensional simple Lie algebras ${\mathcal
M}(m,n)$ of characteristic $5$ depending on two parameters
$m,n\in\mathbb{N}$.

Suppose $\mbox{char}\,F=5$. In \cite{M1,M2}, the algebra ${\mathcal
M}(m,n)$ is described as a graded Lie algebra
$L\,=\,\bigoplus_{i\ge-2}\,L_i$ of dimension $5^{m+n+1}$ whose
graded subalgebra $L_{-2}\oplus L_{-1}$ is isomorphic to a
 five dimensional Heisenberg Lie algebra and $L_0^{(1)}\cong
 W(1;\underline{1})$ as Lie algebras. Moreover,
$L_0=L_0^{(1)}\oplus \mathfrak{z}(L_0)$, $\mathfrak{z}(L_0)=Fz$,
$(\ad z)_{\vert_{L_k}}\,=\,k\cdot\mbox{Id}_{L_k}$ for all
$k\in\mathbb Z$, and $L_{-1}\cong {\mathcal O}(1;\underline{1})/F$
as $W(1;\underline{1})$-modules. It is shown in \cite{M2} that
each Melikian algebra is strongly degenerate and the only
restrictable algebra in the family is ${\mathcal M}(1,1)$ (see
also \cite{Ku''} and \cite{St04} where all derivations of
${\mathcal M}(m,n)$ are determined).

It is stated in \cite{M1} that ${\mathcal M}(1,1)$ is neither a
classical Lie algebra nor a Lie algebra of Cartan type. In
\cite{M2}, Melikian outlines a proof of this statement relying on
properties of $\Z$-gradings in the contact Lie algebra
$K(3;\underline{1})$. An alternative proof will be given below.
Melikian's work showed that the assumption that $p>5$ in the
generalized Kostrikin--Shafarevich conjecture could not be
relaxed.

A few years later Ermolaev observed that $
\mathfrak{g}\,=\,{\mathcal M}(m,n)$ admits a more natural
$\Z$-grading $\mathfrak{g}\,=\,\bigoplus_{i\ge
-3}\,\mathfrak{g}(i)$ that satisfies the conditions (g1), (g2),
(g3) of (2.4) and has the property that $\bigoplus_{i\le
1}\,\mathfrak{g}(i)$, regarded as a local Lie algebra, is
isomorphic to the local Lie algebra associated with a depth $3$
grading of a Lie algebra $\mathcal L$ of type $G_2$. In
particular, the nonpositive part $\bigoplus_{i\le
0}\,\mathfrak{g}(i)$ of $\mathfrak g$ is isomorphic to a maximal
parabolic subalgebra of $\mathcal L$. This observation enabled
Kuznetsov to give in \cite{Ku2} an explicit description of
${\mathcal M}(m,n)$.

Set $\underline{n}:=(m,n)$ and define
$$G_{\bar{0}}\,:=\,\bigoplus_{i\equiv 0\,(\mbox{{\scriptsize
mod}}\,3)} \mathfrak{g}(i),\ \ \ \
G_{\bar{1}}\,:=\,\bigoplus_{i\equiv 1\,(\mbox{{\scriptsize
mod}}\,3)} \mathfrak{g}(i),\ \ \ \
G_{\bar{2}}\,:=\,\bigoplus_{i\equiv 2\,(\mbox{{\scriptsize
mod}}\,3)} \mathfrak{g}(i).$$ Then $
\mathfrak{g}\,=\,G_{\bar{0}}\oplus G_{\bar{1}}\oplus G_{\bar{2}}$
is a $({\Z}/3{\Z})$-grading of $\mathfrak g$. According to
\cite{Ku2}, $$G_{\bar{0}}\oplus G_{\bar{1}}\oplus G_{\bar{2}}\,=\
W(2;\underline{n})\oplus {\mathcal
O}(2;\underline{n})\oplus\widetilde{W}(2;\underline{n})$$ as
vector spaces. Moreover, $G_{\bar{0}}$ is identified with
$W(2;\underline{n})$ as Lie algebras, $G_{\bar{1}}$ is identified
with ${\mathcal O}(2;\underline{n})$ as vector spaces, and
$G_{\bar{2}}$ is identified with
$\widetilde{W}(2;\underline{n})\,=\,\{\widetilde{D}\,\vert\,D\in
W(2;\underline{n})\}$, a vector space copy of
$W(2;\underline{n})$. The Lie product in $\mathfrak g$ is given by
\begin{eqnarray*}
[D,\widetilde{E}] &=&
\widetilde{[D,E]}+2\mbox{div}(D)\,\widetilde{E}, \\\, [D, f] &=&
D(f) - 2\mbox{div}(D)\,f,  \\
{[f_1\widetilde{\partial}_1+f_2\widetilde{\partial}_2,
\,g_1\widetilde{\partial}_1+g_2
\widetilde{\partial}_2]}&=&{f_1g_2-f_2g_1, } \\\, [f,\widetilde{E}]
&=& fE, \\\, [f,g] &=& 2(f\widetilde{{\mathcal
D}}_g-g\widetilde{{\mathcal D}}_f),\ \ \ \ \ {\mathcal
D}_h={\partial}_1(h)\partial_2-\partial_2(h)\partial_1,
\end{eqnarray*}
for all $D,E\in W(2;\underline{n}),\,f,g,h,f_i,g_i\in {\mathcal
O}(2;\underline{n})$. Here $\mbox{div}\,\colon\,W(2;\underline{n})
\rightarrow\, {\mathcal O}(2,\underline{n})$ is the linear map
taking $f_1\partial_1+f_2\partial_2$ to $\partial_1(f_1)+
\partial_2(f_2).$ It follows from the above formulae
that the Lie subalgebra of ${\mathcal M}(m,n)$ generated by the
graded components $\mathfrak{g}(\pm 1)$ is isomorphic to a
classical Lie algebra of type $G_2$.

Assume for a contradiction that ${\mathcal M}(1,1)$ is either
classical or of Cartan type. Since ${\mathcal M}(1,1)$ is strongly
degenerate, simple, and restrictable it must be isomorphic to one
of $W(m;\underline{1})$, $S(m;\underline{1})^{(1)}$,
$H(m;\underline{1})^{(2)}$, $K(m;\underline{1})^{(1)}$. Since
$\dim {\mathcal M}(1,1)=125,$ there is only one option, namely,
${\mathcal M}(1,1)\cong K(3;\underline{1})$. Using the above
multiplication table one can observe that
$\mathfrak{t}_0\,:=\,F(1+x_1)\partial_1\oplus F(1+x_2)\partial_2$
is a torus in ${\mathcal M}(1,1)$ whose centralizer $\mathfrak h$
is a five dimensional Cartan subalgebra of ${\mathcal M}(1,1)$
with the property that
$[\mathfrak{h},[\mathfrak{h},\mathfrak{h}]]=\mathfrak{t}_0$ (see
\cite{P94} for more detail). However, all Cartan subalgebras in
$K(3;\underline{1})$ are abelian, as can be deduced from
\cite{Dem2} and \cite[(7.5)]{St04}. Thus ${\mathcal
M}(1,1)\not\cong K(3;\underline{1})$, and so ${\mathcal M}(1,1)$
is neither classical nor of Cartan type.

Although the Melikian algebras have sporadic nature and can
survive as Lie algebras only at characteristic $5$, they have some
relatives in characteristics $3$ and $2$. This was discovered by
Skryabin \cite{Sk92} and Brown \cite{Br}.

Suppose $\mbox{char}\,F=3$. Each Skryabin algebra $\mathfrak g$ is
equipped with a $\Z$-grading $\mathfrak{g}\,=\,\bigoplus_{i\ge
-4}\,\mathfrak{g}_i$ satisfying the conditions (g1), (g2), (g3) of
(2.4) and one of the three conditions below:
\begin{eqnarray*}
1)\ \ \mathfrak{g}_i=(0)&\mbox{ for } i\le -3&\mbox{ and } \quad\,
\mathfrak{g}_0=\mathfrak{gl}(\mathfrak{g}_{-1}),\\ 2)\ \
\mathfrak{g}_i=(0)&\mbox{ for }  i\le -3&\mbox{ and }\quad\,
\mathfrak{g}_0=\mathfrak{sl}(\mathfrak{g}_{-1}),\\ 3)\ \
\mathfrak{g}_i=(0) &\mbox{ for } i\le -5&\mbox{ and }\quad\,
\mathfrak{g}_{0}=\mathfrak{gl}(\mathfrak{g}_{-1}),\ \
\dim\,\mathfrak{g}_{-4}=3.
\end{eqnarray*}
Moreover, $\dim\,\mathfrak{g}_{-1}=3$ and $\mathfrak{g}_{-2}\cong
\wedge^2\mathfrak{g}_{-1}$ in all cases, and
$\mathfrak{g}_{-3}\cong\wedge^3\mathfrak{g}_{-1}$ in case 3). In
cases 1) and 2), each Skryabin algebra admits a natural
$({\Z}/2{\Z})$-grading $\mathfrak{g}\,=\,G_{\bar{0}}\oplus
G_{\bar{1}}$ such that $G_{\bar{0}}$ is either
$W(3;\underline{n})$ or $S(3;\underline{n};\omega)^{(1)}$ with
$\omega\in {\bf s}(3;\underline{n})$ and $G_{\bar{1}}$ is a nice
irreducible $G_{\bar{0}}$-module. In case 3), each Skryabin
algebra admits a natural $({\Z}/4{\Z})$-grading
$\mathfrak{g}\,=\,G_{\bar{0}}\oplus G_{\bar{1}}\oplus
G_{\bar{2}}\oplus G_{\bar{3}}$ such that
$G_{\bar{0}}\,=\,W(3;\underline{n})$ and each $G_{\bar{i}}$ with
$\bar{i}\ne\bar{0}$ is a nice $G_{\bar{0}}$-module. In all cases,
the Lie bracket in $\mathfrak g$ is given by explicit formulae
involving classical operations with differential forms (see
\cite{Sk92} for more detail).

Now suppose $\mbox{char}\, F=2$. In \cite{Br}, Brown constructed
three series of simple Lie algebras over $F$ one of which relates
closely with the Melikian series.

Following \cite{Br} consider the $({\Z}/3{\Z})$-graded algebra
${\mathcal L}\,=\,{\mathcal L}_{\bar{0}}\oplus{\mathcal
L}_{\bar{1}}\oplus{\mathcal L}_{\bar{2}}$ such that ${\mathcal
L}_{\bar{0}}=W(2;\underline{n})$, ${\mathcal
L}_{\bar{2}}={\mathcal O}(2;\underline{n})$, and ${\mathcal
L}_{\bar{1}}\,=\,\{fu\,\vert\,f\in {\mathcal
O}(2;\underline{n})\}$, a second vector space copy of ${\mathcal
O}(2;\underline{n})$.  The multiplication function
$[\,\cdot\,,\,\cdot\,]\,\colon\,{\mathcal L}\times{\mathcal
L}\rightarrow{\mathcal L}$ satisfies the identity $[x,x]=0$,
agrees with the Lie bracket of $W(2;\underline{n}),$ and has the
following properties:
\begin{eqnarray*}
[D,fu]&=&\mbox{div}(fD)u,\ \ \ \ [D,f]=D(f),\ \ \ \ [fu,gu]=0,\\ \
[fu,g]&=&f{\mathcal D}_{g}, \ \ \ \ \ \ \ \ \ \ \ \
\,[f,g]={\mathcal D}_{g}(f)u
\end{eqnarray*}
(here $f,g\in {\mathcal O}(2;\underline{n})$, $D\in
W(2;\underline{n})$, and ${\mathcal D}_g$ has the same meaning as
before). It is shown in \cite{Br} that $\mathcal L$ is a Lie
algebra carrying a natural $\Z$-grading ${\mathcal
L}\,=\,\bigoplus_{i\ge-4}\,{\mathcal L}_i$ such that ${\mathcal
L}_{-4}\,=\,\mathfrak{z}({\mathcal L})$. The Lie algebra
$\mathfrak{g}\,:=\,({\mathcal L}/\mathfrak{z}({\mathcal
L}))^{(1)}$ is denoted by $G_{2}(2;\underline{n})$. It is simple,
has dimension $2^{|\underline{n}|+2}-2$, and inherits from
$\mathcal L$ a natural $\Z$-grading
$\mathfrak{g}\,=\,\bigoplus_{i\ge -3}\,\mathfrak{g}_i$ satisfying
the conditions (g1), (g2), (g3) of (2.4). Moreover,
$$\mathfrak{g}_0\cong\mathfrak{gl}(\mathfrak{g}_{-1}),\ \ \
\dim\,\mathfrak{g}_{-1}=\dim\,\mathfrak{g}_{-3}=2, \ \ \mbox{ and
}\ \ \mathfrak{g}_{-2}\cong \wedge^2\mathfrak{g}_{-1}.$$ The
$14$-dimensional Lie algebra $G_{2}(2;\underline{1})$ is not
restrictable but can be obtained by reducing modulo $2$  a
nonstandard $\Z$-form of a complex Lie algebra of type $G_2$ (see
\cite{Br} for more detail).

\section{Classification theorems}

One of the main goals of this survey is to announce the following
theorem which, in particular, confirms the original
Kostrikin--Shafarevich conjecture in full generality; see
\cite{PS5}.

\begin{theo}[Classification Theorem]
Let $L$ be a finite dimensional simple Lie algebra over an
algebraically closed field of characteristic $p>3$. Then $L$ is
either a classical Lie algebra or a filtered Lie algebra of Cartan
type or one of the Melikian algebras.
\end{theo}

Our proof of Theorem 7 relies on several earlier classification
results which we are going to formulate. From now on we assume
that $\mbox{char}\,F=p>3$.

The following useful characterization of classical Lie algebras is
due to Seligman and Mills:

\begin{theo}[\cite{MS}]
A Lie algebra $L$ over $F$ is a direct sum of classical simple Lie
algebras if and only if the following conditions hold:
\begin{enumerate}
\item $L$ is perfect and $\mathfrak{z}(L)=(0)$; \item $L$ contains
an abelian Cartan subalgebra $H$ such that
\begin{enumerate}
\item $L=H\oplus\sum_{\alpha\ne 0} L_{\alpha}$ where
$L_{\alpha}\,=\,\{x\in L\,\vert\, [h,x]=\alpha(h)x\ \,
(\forall\,h\in H)\}$; \item if $L_\alpha\neq (0)$, then $\dim\,\,
[L_{\alpha},L_{-\alpha}]=1$; \item if $L_\alpha\neq (0)$ and
$L_\beta\neq (0)$, then $L_{\alpha+k\beta}=(0)$ for some $k\in
\mathbb{F}_p$.
\end{enumerate}
\end{enumerate}
\end{theo}
\noindent
A short proof of the Seligman--Mills theorem based on
the Kac--Moody theory can be found in \cite{Ser}.

The following important theorem allows one to recognize certain
filtered simple Lie algebras:

\begin{theo}[Recognition Theorem] Let $L$ be
a finite dimensional simple Lie algebra over an algebraically
closed field of characteristic $p>3$. Let $$
L=L_{(-s')}\supset\ldots\supset L_{(0)}\supset\ldots\supset
L_{(s)}\supset (0), \ \ \ \ \ \ \ [L_{(i)},L_{(j)}]\subseteq
L_{(i+j)}, $$ be a filtration of $L$ satisfying the following
conditions:
\begin{enumerate}
\item[(a)] $s,s'\ge 1$ and $s'\le s$; \item[(b)] $L_{(0)}/L_{(1)}$
is a direct sum of ideals each of which is either classical simple
or $\mathfrak{gl}(n)$, $\mathfrak{sl}(n)$, $\mathfrak{pgl}(n)$
with $p|n$ or abelian; \item[(c)] $L_{(-1)}/L_{(0)}$ is an
irreducible $L_{(0)}$-module; \item[(d)] for all $j\leq 0$, if
$x\in L_{(j)}$ and $[x,L_{(1)}]\subseteq L_{(j+2)}$, then $x\in
L_{(j+1)}$; \item[(e)] for all $j\geq 0$, if  $x\in L_{(j)}$ and
$[x,L_{(-1)}]\subseteq L_{(j)}$, then $x\in L_{(j+1)}$.
\end{enumerate}
Then $L$ is either classical or is isomorphic as a filtered 
algebra to a Lie algebra of Cartan type or a Melikian algebra
regarded with their natural filtrations.
\end{theo}
\noindent The Recognition Theorem incorporates Wilson's theorem
\cite{Wil76} and earlier results of Kostrikin and Shafarevich
\cite{KS69}. Kac was the first to formulate a version of this
theorem for graded Lie algebras, and he made in \cite{K1} many
deep and important observations towards its proof. One of Kac's
original assumption on the pair $(L_{(-1)},\,L_{(0)})$ was relaxed
by Benkart--Gregory in \cite{BG}. The first complete proof of the
Recognition Theorem for graded Lie algebras was obtained only very
recently by Benkart--Gregory--Premet; see \cite{BGP}. Theorem 9 is
a consequence of this result; see \cite[Section 5]{St04} for more
detail.

Theorems~8 and 9 are fundamental, and most of the classification
proofs rely on them at some stage.

Given a nilpotent Lie subalgebra $H$ of $L$ we denote by
$H_{p}^{\rm tor}$ the unique maximal torus in the $p$-envelope of
$H$ in $\mbox{Der}\,L$. We say that $H$ is {\it triangulable} if
$\mbox{ad}\,h$ is a nilpotent linear operator for any $h\in
H^{(1)}$ (this is the same as to say that $\mbox{ad}\,H$
stabilizes a flag of subspaces in $L$).

We list below a few other classification results which are invoked
frequently. All of them share the assumption that $L$ is a finite
dimensional simple Lie algebra over $F$.

\begin{enumerate}
\item [{\bf 4.1.}] {\it Kaplansky} \cite{Kap}: If $p>3$ and $L$
contains a one dimensional Cartan subalgebra $Ft$ with
$\mbox{ad}\, t$ toral, then $L$ is either $\mathfrak{sl}(2)$ or
$W(1;\underline{1})$.

\smallskip

\item [{\bf 4.2.}] {\it Demushkin} \cite{Dem1, Dem2}, {\it Strade}
\cite[(7.5)]{St04}: If $L$ is a restricted Lie algebra of Cartan
type, then all maximal tori of $L$ have the same dimension and
split into finitely many conjugacy classes under the action of
${\rm Aut}\,L$.

\smallskip

\item [{\bf 4.3.}] {\it Kuznetsov} \cite{Ku1}, {\it Weisfeiler}
\cite{We2}, {\it Skryabin} \cite{Sk2}, {\it Strade} \cite{St04}:
If $p>3$ and $L$ contains a solvable maximal subalgebra, then
either $L\cong\mathfrak{sl}(2)$ or $L\cong W(1;\underline{n})$.

\smallskip

\item [{\bf 4.4.}] {\it Wilson} \cite{Wil77}, {\it Premet}
\cite{P94}: If $H$ is a nontriangulable Cartan subalgebra of $L$,
then $p=5$ and there exist $\mathbb{F}_p$-independent
$\alpha,\beta\in \Gamma(L,H_{p}^{\rm tor})$ and an ideal
$R(\alpha,\beta)$ of the $2$-section $L(\alpha,\beta)$ such that
$$L(\alpha,\beta)/R(\alpha,\beta)\cong {\mathcal M}(1,1).$$

\smallskip

\item [{\bf 4.5.}] {\it Wilson}  \cite{Wil78}, {\it Premet}
\cite{P94}: If $p>3$ and $L$ contains a Cartan subalgebra $H$ with
$\dim\,H_{p}^{\rm tor}=1$, then $L$ is one of $\mathfrak{sl}(2)$,
$W(1;\underline n)$, $H(2;\underline{n};\Phi)^{(2)}$.

\smallskip

\item [{\bf 4.6.}] {\it Block--Wilson} \cite{BW82}, {\it Wilson}
\cite{Wil83}: Suppose $L$ is restrictable and $p>7$. If $L$
contains a toral Cartan subalgebra, then either $L$ is classical
or $L\cong W(n;\underline{1})$.

\smallskip

\item [{\bf 4.7.}] {\it Benkart--Osborn} \cite{BO}: If $L$
contains a one dimensional Cartan subalgebra and $p>7$, then $L$
is either $\mathfrak{sl}(2)$ or $W(1;\underline{n})$ or $L\cong
H(2;\underline{n};\Phi)^{(2)}$ and $\dim L=p^{|\underline{n}|}$.
\end{enumerate}

\noindent The following result of Block--Wilson marked the first
real breakthrough in solving the classification problem for $p>7$.

\begin{theo}[\cite{BW}]
The original Kostrikin--Shafarevich conjecture  is true for $p>7$.
\end{theo}

\noindent Relying heavily on an important intermediate result of
\cite{BW} and the classification techniques of Block--Wilson the
second author was able to generalize Theorem 10, with some support
of R.L.~Wilson (see \cite{St89b, St91, St92, St93, BOSt, St94,
St98a}).
\begin{theo}[Strade 1998] The generalized
Kostrikin--Shafarevich conjecture is true for $p>7$.
\end{theo}

\noindent Large parts of the proof of Theorem 11 go through for
$p>3$ and are incorporated into our proof of Theorem 7.

\section{Principles of the classification}

Let $L$ be a simple Lie algebra over $F$ (recall that
$\mbox{char}\,F=p>3$). As in the characteristic $0$ case we hope
to get more insight into the structure of $L$ by looking at the
root space decomposition of $L$ relative to its Cartan subalgebra
$\mathfrak h$. However, most of the classical results are no
longer valid in our situation. For example, a $(2m+1)$-dimensional
Heisenberg Lie algebra over $F$ admits irreducible representations
 of dimension $p^m$. This implies that Lie's
theorem on solvable Lie algebras fails in characteristic $p$. The
Killing form of any strongly degenerate simple Lie algebra over
$F$ vanishes (see (2.3)). Since all finite dimensional Cartan type
Lie algebras over $F$ are strongly degenerate, Cartan's criterion
is no longer valid in characteristic $p$ either. Cartan
subalgebras of $L$ need not be conjugate under the automorphism
group $\mbox{Aut}\,L$ and, in fact, may have different dimensions
(see our discussion in (2.1)). In characteristic $5$, one can even
expect $L$ to possess nontriangulable Cartan subalgebras (see
({\bf 4.4}) and our discussion in (3.4)).

In general, a nonrestrictable Lie algebra does not possess a
Jordan--Chevalley decomposition. To fix that we embed $L\cong \ad L$
into its semisimple $p$-envelope ${\mathcal L}$ (see ({\bf 2.2.3})).
The Lie algebra ${\mathcal L}\subset\mbox{Der}\,L$ is restricted,
hence admits a Jordan--Chevalley decomposition. By construction,
${\mathcal L}^{(1)}\subseteq L$ (and ${\mathcal L}=L$ if and only if
$L$ is restrictable). We choose a torus $T$ of maximal dimension in
$\mathcal L$ and take a close look the root space decomposition
$$L\,=\,H\oplus\sum_{\alpha\in\Gamma(L,T)}L_{\alpha}$$ of $L$
relative to $T$. Although the subalgebra $H\,=\,\{x\in
L\,\vert\,[t,x]=0\ \ \forall\, t\in T\}$ is nilpotent it is not
always a Cartan subalgebra of $L$ (if $L$ is nonrestrictable, it
may even happen that $H=(0)$).  We wish to gather as much
information as we can on the structure of $1$- and $2$-sections of
$L$ relative to $T$. In characteristic $0$, such information
eventually allows one to determine the global structure of $L$.

In characteristic $p$, the local analysis is much more involved.
There are a number of reasons for that. To mention just a few, the
$1$-sections of $L$ relative to $T$ are no longer ``reductive'' and
their irreducible representations are hard to describe. Some tori of
maximal dimension in $\mathcal L$ are unsuitable for our purposes,
and a lot of effort is spent on optimizing a randomly chosen $T$ by
using generalized Winter exponentials; see (2.1). In the course of
the proof one has to make various sophisticated choices of maximal
subalgebras, carry out detailed computations in Lie algebras of
small rank, and study central extensions of such algebras and their
irreducible representations.

For any $\alpha\in\Gamma(L,T)$ the semisimple quotient $L[\alpha]$
of the 1-section $L(\alpha)$ is either zero or $\mathfrak{sl}(2)$
or $W(1;\underline{1})$ or the inclusion
$$H(2;\underline{1})^{(2)}\subset L[\alpha]\subset
H(2;\underline{1})$$ holds (this follows from ({\bf 4.5})).
Accordingly we call $\alpha$ {\it solvable, classical, Witt} or
{\it Hamiltonian}. It is not difficult to show that the radical of
$L(\alpha)$ is $T$-stable. So $T$ acts as derivations on
$L[\alpha]$ and $L[\alpha]^{(2)}$. Following Block--Wilson we say
that $\alpha$ is a {\it proper} root if either
$L[\alpha]\in\{(0),\mathfrak{sl}(2)\}$ or $L[\alpha]$ is of Cartan
type and the standard maximal subalgebra of $L[\alpha]^{(2)}$ is
$T$-invariant. If $\alpha$ is not a proper root we say that
$\alpha$ is {\it improper}.

The main intermediate result of \cite{BW} is a classification of
all simple Lie algebras of absolute toral rank $2$ for $p>7$.
Combining this classification with a version of {(\bf 4.4}) for
$p>7,$ Block and Wilson succeeded to describe the semisimple
quotients of all 2-sections in a restricted simple Lie algebra.
Having achieved that they proceed as follows:

The description of the quotients $L[\alpha]$ mentioned above
implies that each $1$-section $L(\alpha)$ contains a unique
subalgebra $Q(\alpha)$ with $H\subset Q(\alpha)$ and
$\dim\,Q(\alpha)=\dim\,L(\alpha)-e(\alpha)$, where $$e(\alpha)=
\left\{\begin{array}{lll} 0 & \mbox{if $\alpha$ is solvable or
classical},\\ 1 & \mbox{if $\alpha$ is Witt},\\ 2 & \mbox{if
$\alpha$ is Hamiltonian}.
\end{array} \right.$$ The subalgebra
$Q(\alpha)$ is solvable if $\alpha$ is solvable or Witt, and
$Q(\alpha)/\mbox{rad}\,Q(\alpha)\cong\mathfrak{sl}(2)$ if $\alpha$
is classical or Hamiltonian. In all cases, $Q(\alpha)$ is
$T$-invariant if and only if $\alpha$ is a proper root of $L$.
Generalized Winter exponentials are now used to ``optimize'' $T$.
A torus $T\subset\mathcal L$ is called {\it optimal} if $\dim\,
T=MT({\mathcal L})$ and the number of proper roots in
$\Gamma(L,T)$ is maximal possible. Using their description of the
semisimple quotients $L(\alpha,\beta)/\rad L(\alpha,\beta)$ Block
and Wilson prove that in the restricted case {\it all} roots of
$L$ relative to an optimal torus $T\subset\mathcal L$ are proper.
They then look again at the 2-sections of $L$ relative to $T$ to
prove that the $T$-invariant subspace
$$Q=Q(L,T)\,:=\,\sum_{\alpha\in\Gamma(L,T)} Q(\alpha)$$ is a Lie
subalgebra of $L$. The rest of the proof is straightforward. If
$Q=L$, Block and Wilson show that the Seligman--Mills theorem
applies to $L$. So $L$ is classical in this case. If $Q\ne L$,
they show that $Q$ can be embedded into a maximal subalgebra
satisfying the conditions of the Recognition Theorem.

For an arbitrary simple $L,$ the second author used the
Block--Wilson classification of simple Lie algebras of rank $2$ to
obtain a list of all possible $T$-semisimple quotients of the
$2$-sections of $L$ (this list is longer than in the restricted
case). He then succeeded to optimize $T$ in $\mathcal L$ and in the
joint work with Benkart and Osborn \cite{BOSt} constructed a large
Lie subalgebra $Q=Q(L,T)$ of $L$. However, the final parts of the
proof in the general case are {\it much} more involved; see
\cite{St91, St93, St94, St98a}. Essentially, this is due to the fact
that $\mathcal L$ is no longer simple. Since optimal tori may lie
outside $L$ some 3-sections have to be thoroughly investigated.

It turned out that if all regular Cartan subalgebras of $\mathcal
L$ are triangulable, then the final parts of the second author's
classification go through for $p>3$ after a proper modification.
This modification is carried out in \cite{PS4, PS5}, thus settling
the remaining case $p=7$ of the generalized Kostrikin-Shafarevich
conjecture. If $\mathcal L$ contains a nontriangulable regular
Cartan subalgebra, then \cite[Theorem A] {PS4} and ({\bf 4.4})
imply that $p=5$ and one of the semisimple quotients
$L(\alpha,\beta)/\rad L(\alpha,\beta)$ is isomorphic to ${\mathcal
M}(1,1)$. This situation is investigated in \cite{PS6}, the last
paper of the series. The main result of \cite{PS6} states that $L$
is then isomorphic to a Melikian algebra ${\mathcal M}(m,n)$.

The hardest part of our proof of Theorem 7 is the classification
of the simple Lie algebras of absolute toral rank $2$ and the
description of the $2$-sections of $L$ relative to $T$. The former
is obtained in \cite{PS1, PS2, PS3} while the latter is carried
out in \cite{PS4}.

\medskip

Below we outline our arguments in the rank $2$ case.

\medskip

\noindent (A) From now on we assume that $L$ is a nonclassical
simple Lie algebra of absolute toral rank $2$ and $\mathcal L$ is
the semisimple $p$-envelope of $L$; see ({\bf 2.2.3}) and
Definition~3 in (2.2). In view of ({\bf 4.5}) we may assume that
for any maximal torus $T\subset {\mathcal L}$ the centralizer
$H=\mathfrak{c}_{L}(T)$ has the property that $\dim\,H_{p}^{\rm
tor}=2$ (in particular, it can be assumed that all maximal tori in
$\mathcal L$ are two dimensional). Finally, we may assume that all
simple Lie algebras $\mathfrak g$ with $TR(\mathfrak{g})=2$ and
$\dim\,\mathfrak{g}<\dim\, L$ are known. Our ultimate goal is to
prove that $L$ admits a filtration satisfying the conditions of
the Recognition Theorem. However, at the beginning of the
investigation {\it any} long filtration invariant under a two
dimensional torus in $\mathcal L$ would do. Thus we have to
address the following

\smallskip

\noindent {\bf Problem.} {\it Find a {\rm long} standard filtration
in $L$ stable under the action of a maximal torus in $\mathcal L$}.

\smallskip

\noindent This problem is solved in \cite{PS1} by producing a {\it
root sandwich} in $L$, that is a nonzero sandwich element $c\in L$
such that $[T,c] \subset Fc$ for some torus $T$ of maximal
dimension in $\mathcal L$. The set of all such sandwiches is
denoted by ${\mathcal S}(L,T)$. Adopting the method used in
\cite{P86a, P86c} for proving Kostrikin's conjecture we first show
that under some mild assumptions on a $1$-section of $L$ there
exists a nonzero
$$x\in H\cup \bigcup_{\gamma\in\Gamma(L,T)}L_\gamma$$ such that
$(\ad x)^3=0$. Then we use some techniques from \cite{Be, K69} and
the theory of finite dimensional Jordan algebras to find a root
sandwich $c$. More precisely, we prove

\begin{theo}[{\cite{PS1}}]
Let $\mathfrak g$ be a simple Lie algebra of absolute toral rank
$2$ over $F$. Then either $\mathfrak g$ is classical or
$\mathfrak{g}\cong H(2;\underline{1}; \Phi)^{(2)}$ with $\dim
\mathfrak{g}=p^2-1$ or there exists a two dimensional torus
$\mathfrak t$ in the semisimple $p$-envelope of $\mathfrak g$ such
that ${\mathcal S}(\mathfrak{g},\mathfrak{t})\ne \emptyset$.
\end{theo}

Having found a root sandwich $c\in L$ we now observe that any
maximal subalgebra $L_{(0)}$ of $L$ containing $H+
\mathfrak{c}_{L}(c)$ gives rise to a long $T$-invariant filtration
of $L$. Indeed, let $L_{(-1)}$ be any $L_{(0)}$-stable subspace of
$L$ such that $L_{(-1)}\supsetneq L_{(0)}$ and $L_{(-1)}/L_{(0)}$
is an irreducible $L_{(0)}$-module. The $L_{(0)}$-module
$L_{(-1)}$ is $H$-stable, hence $T$-stable (for $T=H_{p}^{\rm
tor}$). Therefore, so are all components of the standard
filtration associated with the pair $(L_{(-1)},\,L_{(0)})$; see
(2.4) for more detail. Since $[L_{(-1)},c]\subset
\mathfrak{c}_{L}(c)\subset L_{(0)}$ we have that $0\ne c\in
L_{(1)}$.

\medskip

\noindent (B) Next we investigate the graded Lie algebra $G:=\gr
L$. Let $M(G)$ denote the largest ideal of $G$ contained in
$\sum_{i<-1}\, G_i$, and $\bar G:=G/M(G)$. By a theorem of
Weisfeiler \cite{We1}, the Lie algebra $\bar G$ is semisimple and
has a unique minimal ideal, denoted $A(\bar G)$. Furthermore,
$\bar G$ inherits a natural grading from $G$ which satisfies the
conditions (g1) - (g4). Note that finite dimensional semisimple
Lie algebras  over $F$ need not be direct sums of simple ideals
(in fact, {\it simple} ideals may not exist at all). The structure
of semisimple modular Lie algebras was determined by Block in
\cite{Bl}. The following important theorem describes the structure
of a semisimple Lie algebra with a unique minimal ideal:

\begin{theo}[Block's Theorem] Let $\mathfrak g$ be a finite dimensional
semisimple Lie algebra over an algebraically closed field of
characteristic $p>0$ and suppose that $\mathfrak g$ contains a
unique minimal ideal, $I$ say. Then there exist an
$r\in\mathbb{N}_0$ and a simple Lie algebra $\mathfrak s$ such
that $ I\cong \mathfrak{s}\otimes {\mathcal O}(r;\underline 1)$ as
Lie algebras. Moreover, $\mathfrak{g}\cong {\rm
ad}_I\,\mathfrak{g}$ and
$$(\ad\mathfrak{s})\otimes {\mathcal O}(r;\underline{1})\subset {\rm
ad}_I\,\mathfrak{g}\subset (\Der \mathfrak{s})\otimes {\mathcal
O}(r;\underline 1) \rtimes {\rm Id}_{\mathfrak s}\otimes
W(r;\underline 1).$$
\end{theo}
In a sense, the above-mentioned theorem of Weisfeiler can be
regarded as a graded version of Block's theorem; see
\cite[(3.5)]{St04} for more detail. In \cite{PS2}, we show that
our maximal subalgebra $L_{(0)}$ can be chosen such that
$$\mbox{either }\ G_2\ne (0)\ \,\mbox{   or }\ \,
[[G_{-1},G_1],G_1]\ne(0).$$ In this case, Weisfeiler's theorem
says that $A(\bar G)\,=\,\bigoplus_{i}A(\bar G)_i$ where $A(\bar
G)_i\,=\,A(\bar G)\cap\bar G_i$ and there exist a graded simple
Lie algebra $S\,=\,\bigoplus_{i}S_i$ and an integer $m\ge 0$ such
that
$$A(\bar G)\cong S\otimes {\mathcal
O}(m;\underline{1})\,=\,\bigoplus_{i\in\Z}\, (S_i\otimes {\mathcal
O}(m;\underline{1}))$$ as graded Lie algebras. In \cite{PS2} we
show that $m\le 1$. Moreover, we prove that if $m=1$, then the
absolute toral rank of $S$ drops. In view of ({\bf 2.2.4}) and
({\bf 4.5}) the equality $m=1$ implies that $S$ is one of
$\mathfrak{sl}(2)$, $W(1;\underline 1)$, $H(2;\underline
1)^{(2)}$.

We first consider the case where $m=1$. By Block's theorem, we
then have an embedding
\begin{eqnarray*}
\bar G&\hookrightarrow &(\Der S)\otimes {\mathcal O}(1;\underline
1)\rtimes{\rm Id}_S\otimes  W(1;\underline 1).
\end{eqnarray*}
In view of a conjugacy theorem proved in \cite{PS2} along comes an
induced embedding of tori $$ T\hookrightarrow T_0\otimes 1+{\rm
Id}_S\otimes Fz\partial,\quad z\in\{x,1+x\}, $$ where $T_0$ is a
one dimensional torus in $\Der S$. We then show that $T$ and
$L_{(0)}$ can be chosen such that $S\cong H(2;\underline 1)^{(2)}$
as graded Lie algebras, where $H(2;\underline{1})^{(2)}$ is
regarded with its grading of type $\underline{1}$; see (3.2). In
particular, $S_0\cong \mathfrak{sl}(2)$ and $S_{-k}=(0)$ for $k\ge
2$. We also show that $M(G)=(0)$. This information enables us to
conclude, eventually, that $p=5$ and $L\cong{\mathcal M}(1,1)$.

\medskip

\noindent (C) From now on we may assume that $m=0$. Using the
inequality $TR(G)\le TR(L)$, proved in \cite{Sk3}, we  show that
$TR(S)= 2$. We now wonder whether $S$ is listed in the
Classification Theorem.

First we observe that the root sandwich $c\in L_{(1)}$ gives rise
to a nonzero sandwich element of $G$ contained in the graded
component $G_l$ for some $l\ge 1$. Since
$M(G)\,\subset\,\bigoplus_{i<-1}\,G_i$, it follows that the Lie
algebra $\bar G$ must be strongly degenerate. Since $S\subset\bar
G\subset \Der S$, it follows that $S$ is not a classical Lie
algebra.

Next we observe that the quotient space $\bar M:=M(G)/M(G)^2$ is a
$\bar G$-module, hence an $S$-module. Let $S_p$ denote the
$p$-envelope of $S$ in $\Der\, S$. We show in \cite{PS3} that any
composition factor $V$ of the $S$-module $\bar M$ can be viewed in a
natural way as a restricted $S_p$-module and $T$ can be identified
with a two dimensional torus in $S_p$. Since $H\subset L_{(0)}$ it
must be that $$0\not\in\Gamma^w(V,T).$$ On the other hand, we show
in \cite{PS3} that if $S$ is isomorphic to one of
$S(3;\underline{1})^{(1)}$, $H(4;\underline{1})^{(2)}$,
$K(3;\underline{1})$, ${\mathcal M}(1,1)$, $H(2;(2,1))^{(2)}$, then
$T$ has weight $0$ on {\it any} finite dimensional restricted
$S_p$-module. This implies that $M(G)=(0)$ if $S$ is one of these
Lie algebras. A slight modification of the argument shows that
$M(G)=(0)$ if $S$ is one of $W(2;\underline{1})$, $W(1;2)$,
$H(2;\underline{1};\Phi)^{(2)}$. Since $TR(S)=2$ we deduce the
following: $$\mbox{if $S$ is known, then } M(G)=(0).$$

\medskip

Suppose $S$ is known. Then $\bar G\cong G=\gr L$, hence $L$ is a
filtered deformation of $\bar G\subset \Der S$. So there exists a
Lie algebra $\mathfrak L$ over the polynomial ring $F[t]$ such
that
$$ \mathfrak{L}/(t-\lambda)\mathfrak{L}\,\cong L \ \mbox{ if }\,
\lambda\ne0,\, \mbox{ and }\ \mathfrak{L}/t\mathfrak{L}\,\cong\,
{\bar G}\supset S. $$

\medskip

Suppose $S$ is a Melikian algebra. Since $TR(S)=2$ we then have
$S={\mathcal M}(1,1)$. By \cite{Ku''}, all derivations of
${\mathcal M}(1,1)$ are inner (see also \cite[(7.1)]{St04}). So it
must be that $G\cong S$. We already mentioned in (3.4) that
${\mathcal M}(1,1)$ contains a two dimensional torus $
\mathfrak{t}_0$ whose centralizer $\mathfrak h$ is a
nontriangulable Cartan subalgebra of ${\mathcal M}(1,1)$. As
$TR(S)=2$, the Cartan subalgebra $\mathfrak h$ is regular in $S$.
As all regular Cartan subalgebras of a finite dimensional
restricted Lie algebra have the same dimension (see (2.1)) we can
lift $\mathfrak h$ to a nontriangulable Cartan subalgebra of
minimal dimension in $\mathfrak{L}\otimes_{F[t]} F(t)$. We then
use a deformation argument to show that $L$ contains a
nontriangulable Cartan subalgebra as well. Using ({\bf 4.4}) we
finally conclude that $L\cong {\mathcal M}(1,1)$.

Suppose $S\cong X(m;\underline{n})^{(2)}$ where $X\in\{W,S,H,K\}$.
Any grading of a Lie algebra $\mathfrak g$ is induced by the
action of a one dimensional torus of the algebraic group
$\mbox{Aut}\, \mathfrak g$. Each such torus is contained in a
maximal torus of $\mbox{Aut}\,\mathfrak g$. The conjugacy theorem
for maximal tori of algebraic groups enables us to prove that any
grading of $S$ is obtained by assigning certain integral weights
to the elements of a generating set of the divided power algebra
${\mathcal O}(m;\underline{n})$. This procedure also describes the
gradings of $\Der S$ and provides valuable information on gradings
of $G$ (for $G$ can be regarded as a graded subalgebra of $\Der
S$). It turns out that very few gradings of $G$ can satisfy the
conditions (g1), (g2), (g3). Taking graded Cartan type Lie
algebras of rank $2$ one at a time we show that our choice of
$L_{(0)}$ (and $T$) forces the grading of $S\cong
X(m;\underline{n})^{(2)}$ to be standard. At this point Wilson's
theorem enables us to conclude that $L$ is a filtered Lie algebra
of Cartan type.

If $S$ is a filtered Cartan type Lie algebra not considered
before, that is one of type $H(2;\underline{1};\Phi)^{(2)}$, then
$S$ is nonrestrictable of dimension $p^2-1$ or $p^2$. In this
case, $S_p$ is known to possess a two dimensional toral Cartan
subalgebra $\mathfrak t$ with the property that $\dim\,S_\gamma=1$
for all $\gamma\in\Gamma(S,\mathfrak{t})$. This information and an
intermediate result of \cite{BW82} (applicable for $p>3$ in view
of ({\bf 4.5})) allow us to show that $L$ too is of Cartan type.

\medskip

\noindent (D) It remains to consider the case where $S\cong
G=\mbox{gr}\,L$ is a minimal counterexample to our theorem. At
this stage we may also assume that passing from $L$ to $G$ {\it
always} produces unknown simple graded Lie algebras (subject to
certain conditions on $T$ and $L_{(0)}$). We use this as a
technical tool for improving $L_{(0)}$ and obtaining more
information on the structure of $\Gamma\big(\sum_{i<0}\,G_i,
T\big)$. Given $\alpha\in\Gamma(G,T)$ we set $K_{\alpha}\,
:=\,\{x\in G_\alpha\,\vert\,\alpha([x,G_{-\alpha}])=0\}$ and
denote by $K'(G,T,\alpha)$ the Lie subalgebra of $G$ generated by
all $K_{i\alpha}$ with $i\in\mathbb{F}_p^*$. It follows from the
main result of \cite{PS2} that $K'(G,T,\alpha)$ is a triangulable
subalgebra of $G$.

The most important task for us now is to determine the graded
component $G_0$. From \cite{Sk2} we know that the radical of $G_0$
is abelian, while ({\bf 4.3}) entails that $G_0$ is nonsolvable.
Thus if $\mbox{rad}\,G_0\ne (0)$, then $\bar G_0\,:=\,
G_0/\mbox{rad}\,G_0$ has absolute toral rank $1$. Moreover,  it
follows from ({\bf 4.5}) that $\bar G_0$ is either $
\mathfrak{sl}(2)$ or $W(1;\underline{1})$ or the inclusion
$H(2;\underline{1})^{(2)}\subset \bar G_0\subset
H(2;\underline{1})$ holds. Combining some representation theory
with the fact that $K'(G,T,\alpha)$ is triangulable (see
\cite{PS2}), we show after a detailed analysis that either
$G_0\cong W(1;\underline 1)\ltimes {\mathcal O}(1;\underline 1)$
(a natural semidirect product) or the radical of $G_0$ is one
dimensional and central, and the extension $$0\to \rad G_0\to
G_0\to \bar G_0\to 0$$ splits. If $G_0$ is semisimple with a
unique minimal ideal $I$, then Block's theorem says that $I\cong
\mathfrak{s}\otimes {\mathcal O}(r;\underline{1})$ for some simple
Lie algebra $\mathfrak s$. If $r>0$ we prove that $\mathfrak s$
has absolute toral rank $1$ and there are a vector space $V$ over
$F$ and a linear isomorphism
$$G_{-1}\,\stackrel{\sim}{\longrightarrow} \,V\otimes
{\mathcal O}(k;\underline l)$$ such that ${\mathcal
O}(k,\underline l)\cong {\mathcal O}(r;\underline 1)$ as algebras
and the action of $G_0$ on $G_{-1}$ is induced by a Lie algebra
embedding
$$G_0\hookrightarrow\,\mathfrak{gl}(V)\otimes {\mathcal O}(k;\underline l)\rtimes\,{\rm
Id}_V \otimes W(k;\underline l). $$ Moreover, $\pi(G_0)$, the
image of $G_0$ under the canonical projection
$$\pi\,\colon\,\mathfrak{gl}(V)\otimes {\mathcal O}(k;\underline l)\rtimes\,{\rm
Id}_V\otimes W(k;\underline l)\,\longrightarrow\,W(k;\underline
l),$$ is {\it transitive}, that is has the property that
$\pi(G_0)+W(k;\underline l)_{(0)}\,=\,W(k;\underline l)$. Using
the simplicity of $G$ and Cartan prolongation techniques inspired
by earlier work of Kuznetsov (see e.g. \cite{Ku1}) we show that
$\pi(G_0)$ is an ${\mathcal O}(k;\underline l)$-submodule of
$W(k;\underline{l})$. The transitivity of $\pi(G_0)$ now forces
$\pi(G_0)=W(k;\underline l)$, while toral rank considerations
yield $k=1$, $\underline{l}=\underline{1}$. This enables us to
prove that $$G_0\cong \mathfrak{s}\otimes {\mathcal
O}(1;\underline 1)\rtimes{\rm Id}_{\mathfrak s}\otimes
W(1;\underline 1),$$ where $\mathfrak s$ is either
$\mathfrak{sl}(2)$ or $W(1;\underline 1)$. As a consequence, we
obtain that $G_0$ belongs to a short list of known linear Lie
algebras.

Considering algebras from this list one at a time we show that $T$
can be chosen such that all roots in $\Gamma(G,T)$ are proper.
This allows us to obtain much better estimates for $\dim\,
G_{i,\gamma}$ with $i<0$ and $\gamma\in\Gamma(G,T)$. We use this
new information to show that either $G_0$ is a classical Lie
algebra of rank $2$ or the $p$-envelope ${\mathcal G}_0$ of $G_0$
in $\Der G$ is isomorphic to $\mathfrak{gl}(2)$ as restricted Lie
algebras.

Let $G'$ denote the Lie subalgebra of $G$ generated by $G_{\pm 1}$
and $M(G')$ the maximal graded ideal of $G'$ contained in
$\sum_{i<0}\,G_i$. If $M(G')\ne (0)$ we combine the Recognition
Theorem with some representation theory of Cartan type Lie
algebras to show that ${\mathcal G}_0\cong \mathfrak{gl}(2)$ and
$G'/M(G')$ is classical of type $A_2$, $C_2$ or $G_2$. We then use
the representation theory of algebraic groups to show that this
cannot happen. As a result, $M(G')=(0)$. Then the Recognition
Theorem applies to $G$ itself, showing that $G$ is known. This
contradiction proves the Classification Theorem in the rank $2$
case.

\section{Some open problems}

The classification problem in characteristics $2$ and $3$ is wide
open. Since our knowledge of finite dimensional simple Lie
algebras over algebraically closed fields of characteristics $2$
and $3$ is very limited, it is not clear at present whether a
complete classification of such algebras can ever be achieved. As
indicated in our discussion at the end of (3.3) the classification
of Hamiltonian forms in ${\bf h}_1(m,\underline{n})$ was reduced
by Skryabin to a certain problem of linear algebra. Luckily, the
problem turned out to be tame. But if it turned out to be wild, we
would never have a {\it complete} classification in characteristic
$p>3$.

\smallskip

The first three items will address issues in characteristics $2$ and
$3$.

\smallskip

\noindent {\bf Conjecture~1}. {\it The automorphism group of any
finite dimensional simple Lie algebra over an algebraically closed
field of characteristic $p>0$ is infinite.}

\smallskip

\noindent For $p>3$, one can easily deduce Conjecture from the
results in \cite{P86a, P86c} or, alternatively, from Theorem~7.
However, the conjecture remains wide open for $p\in\{2, 3\}$.

\smallskip

In \cite{Sk3}, Skryabin proved that any finite dimensional simple
Lie algebra of absolute toral rank one over an algebraically
closed field of characteristic $3$ is isomorphic to either
$\mathfrak{sl}(2)\cong W(1;\underline{1})$ or
$\mathfrak{psl}(3)\cong H(2;\underline{1})^{(2)}$. He also proved
in \textit{loc. cit.} that no finite dimensional simple Lie
algebras of absolute toral rank $1$ exist in characteristic $2$.

\smallskip

\noindent {\bf Problem~1}. {\it Classify all finite dimensional
simple Lie algebras of absolute toral rank two over algebraically
closed fields of characteristics $2$ and $3$.}

\smallskip
In characteristic $2$, strong results closely related to Problem~1
are obtained by A. Grishkov and the first author (work in
progress). We are unaware of any ongoing work on the
characteristic $3$ case of Problem~1.

\smallskip

As mentioned in \cite{BGP}, it would be very useful to have a
version of the Recognition Theorem for graded Lie algebras of
characteristics $2$ and $3$.

\smallskip
\noindent {\bf Problem~2}. {\it Classify all finite dimensional
graded Lie algebras $G=\bigoplus_{i\in\Z}\,G_i$ over algebraically
closed fields of characteristics $2$ and $3$ that satisfy the
conditions (g1) - (g4) of (2.4) and have the property that $G_0$
is isomorphic to the Lie algebra of a reductive group.}

\smallskip

The last four items will deal, mainly, with the case where $p>3$.

\smallskip

\noindent {\bf Problem~3}. {\it Determine the absolute toral rank of
all finite dimensional simple Lie algebras over algebraically closed
fields of characteristic $p>3$.}

\smallskip
 One should stress here that the value of $TR(L)$ is known for many simple Lie
algebras $L$. In particular, it is known for all restricted Lie
algebras of Cartan type. Some results related to Problem~3 can be
found in \cite{BKK}. The most interesting open case of Problem~3
is the case where $L=H(m;\underline{n};\omega_{A,B})^{(2)}$ and
$\omega_{A,B}\in {\bf h}_1(m;\underline{n})$ is such that $\det
B=0$.

\smallskip

\noindent {\bf Problem~4}. {\it Determine the automorphism groups
of all finite dimensional simple Lie algebras over algebraically
closed fields of characteristic $p>3$. In particular, is it true
that any finite dimensional simple Lie algebra $L$ admits a
$\Z$-grading $L=\bigoplus_{i\in\Z}\,L_i$ with $L_0\ne L$?
Equivalently, is it true that the connected component of the
algebraic group ${\rm Aut}\,L$ is not unipotent?}

\smallskip
There are many examples of simple Lie algebras with {\it solvable}
automorphism groups; in fact, such algebras occur in all four
Cartan series. Probably, the most interesting open case of
Problem~4 is the case where
$L=H(m;\underline{n};\omega_{A,B})^{(2)}$ and $\omega_{A,B}\in
{\bf h}_1(m;\underline{n})$ is such that $\det B\ne 0$.

\smallskip

The next problem was suggested to the first author by R. Guralnick.

\smallskip

\noindent {\bf Question~1} (cf. \cite[Question~2.3]{GKPS}). {\it
Is it true that any finite dimensional simple Lie algebra over an
algebraically closed field of characteristic $p>0$ can be
generated by two elements?}

\smallskip

For $p\in\{2,3\}$ Question~1 is out of reach at the moment. However,
for $p>3$, finite dimensional simple Lie algebras are likely to
enjoy a much stronger property which is nowadays referred to as
``one and a half generation''.

\smallskip
\noindent {\bf Problem~5} (cf. \cite[Question~2.4]{GKPS}). {\it
Let $L$ be a finite dimensional simple Lie algebra over an
algebraically closed field of characteristic $p>3$. Use Theorem~7
to prove that for any nonzero $x\in L$ there is $y\in L$ such that
$L=\langle x,y\rangle$.}

\smallskip

The simplicity assumption on $L$ in Problem~5 is crucial. Indeed,
analyzing the semidirect products
$${\mathcal L}(\mathfrak{g},m):=\,\big(\text{Id}_{\mathfrak g}\,\otimes
{\mathcal D}\big)\ltimes \big(\mathfrak{g}\otimes{\mathcal
O}(m;\underline{1})\big),$$ where $\mathfrak g$ is a finite
dimensional simple Lie algebra over $F$ and $\mathcal D$ is the
commutative subalgebra of $W(m;\underline{1})$ spanned by
$\partial_1,\ldots,\partial_m$, one can observe that for any
natural number $n$ there exists a finite dimensional {\it
semisimple} Lie algebra $L$ over $F$ such that the set
$$S_n(L):=\,\{x\in L\,|\,\,\langle x,y_1,\ldots,y_n\rangle\, \mbox{ is
solvable for all }\  y_1,\ldots, y_n\in L\}$$ is nonzero. More
precisely, it is not hard to see that for the semisimple Lie
algebra $L={\mathcal L}(\mathfrak{g}, n+1)$ one has
$$\mathfrak{g}\otimes x_1^{p-1}\cdots x_n^{p-1}x_{n+1}^{p-1}\subseteq
S_n(L).$$ This is in sharp contrast with the situation for finite
groups; see \cite[Theorem~1.1]{GKPS} for more detail.

\bibliographystyle{amsalpha}

 \end{document}